\documentclass{amsart}
\usepackage[mathscr]{eucal}
\usepackage{ulem}
\usepackage[utf8]{inputenc}
\usepackage{bbm}
\usepackage{dsfont}
\usepackage[reversemp,
    paperwidth=210mm,
    paperheight=297mm,
    top={26mm},
    headheight={5.5pt},
    headsep={5.6mm},
    text={150mm,237mm},
    marginparsep=5mm,
    marginparwidth=12mm,
    bindingoffset=10mm,
    footskip=10mm]{geometry}
\usepackage{color}
\usepackage{mathrsfs}
\usepackage{amssymb,amsmath}
\usepackage{calligra}
\usepackage{fontenc}
\usepackage{amsbsy}
\numberwithin{equation}{section}
\DeclareMathAlphabet{\mathpzc}{OT1}{pzc}{m}{it}
\newcommand{\W}{{\normalfont\textsf{W}^{2}_\sharp}}
\newcommand{{\M}}{{\normalfont\textsf{L}^{2}_\sharp}}

\newcommand{\bsfW}{\boldsymbol{\mathsf W}}
\newcommand{\bsfV}{\boldsymbol{\mathsf V}}
\newcommand{\bsfU}{\boldsymbol{\mathsf U}}

\newcommand{\bsfu}{\boldsymbol{\mathsf u}}

\newcommand{\sub}[1]{{\mbox{\footnotesize $#1$}}}
\usepackage{graphicx}
\graphicspath{%
    {converted_graphics/}% inserted by PCTeX
    {/}% inserted by PCTeX
}

\newcommand{\sfp}{{\sf p}}

\newcommand{\Cdot}{\!\cdot\!}

\newcommand{\Div}{\mbox{\rm div}\,}

\newcommand{\supp}{\mbox{\rm supp}\,}
\newcommand{\curl}{\mbox{\rm curl}\,}

\newcommand{\Int}[2]{{\displaystyle \int_{ #1}^{ #2}}}
\newcommand{\Lim}[1]{{\displaystyle \lim_{ #1}}}

\newcommand{\Sup}[1]{{\displaystyle \sup_{#1}}}

\newcommand{\Sum}[2]{{\displaystyle \sum_{#1}^{#2}}}

\newcommand{\Frac}[2]{\displaystyle{\frac{\displaystyle{#1}}{\displaystyle{#2}}}}

\newcommand{\beea}{\begin{eqnarray}}
\newcommand{\eeea}{\end{eqnarray}}

\newcommand{\bfe}{{\mbox{\boldmath $e$}} }

\newcommand{\bfz}{{\mbox{\boldmath $z$}} }
\newcommand{\0}{{\mbox{\boldmath $0$}} }

\newcommand{\BF}{\begin{footnotesize}}
\newcommand{\EF}{\end{footnotesize}}
\newcommand{\ode}[2]{{\displaystyle \frac{\mbox{$d #1$}}{\mbox{$d #2$}}}}

\newcommand{\bi}{\begin{itemize}}
\newcommand{\ei}{\end{itemize}}
\newcommand{\ed}{\end{document}}
\newcommand{\be}{\begin{equation}}
\newcommand{\bestar}{\begin{equation*}}
\newcommand{\ba}{\begin{array}}
\newcommand{\ea}{\end{array}}
\newcommand{\eestar}{\end{equation*}}
\newcommand{\eeq}[1]{\label{eq:#1}\end{equation}}
\newcommand{\real}{\mathbb R}

\newcommand{\nat}{{\mathbb N}}
\newcommand{\bfpsi}{\mbox{\boldmath $\psi$}}

\newcommand{\bfxi}{\mbox{\boldmath $\xi$}}

\newcommand{\bfphi}{\mbox{\boldmath $\varphi$}}

\newcommand{\bfv}{{\mbox{\boldmath $v$}} }
\newcommand{\bfu}{{\mbox{\boldmath $u$}} }

\newcommand{\bfw}{{\mbox{\boldmath $w$}} }
\newcommand{\bff}{{\mbox{\boldmath $f$}} }
\newcommand{\bfa}{{\mbox{\boldmath $a$}} }

\newcommand{\bfc}{{\mbox{\boldmath $c$}} }

\newcommand{\bfS}{{\mbox{\boldmath $S$}} }

\newcommand{\bfM}{{\mbox{\boldmath $M$}} }

\newcommand{\bfh}{{\mbox{\boldmath $h$}} }

\newcommand{\calb}{{\mathcal B}}
\newcommand{\calc}{{\mathcal C}}
\newcommand{\cald}{{\mathcal D}}

\newcommand{\calf}{{\mathcal F}}
\newcommand{\calg}{{\mathcal G}}
\newcommand{\calh}{{\mathcal H}}
\newcommand{\cali}{{\mathcal I}}

\newcommand{\calk}{{\mathcal K}}
\newcommand{\call}{{\mathcal L}}

\newcommand{\calx}{{\mathcal X}}
\newcommand{\caly}{{\mathcal Y}}

\newcommand{\calz}{{\mathcal Z}}
\newcommand{\bfsigma}{\mbox{\boldmath $\sigma$}}

\newcommand{\bftau}{\mbox{\boldmath $\tau$}}

\newcommand{\bfV}{{\mbox{\boldmath $V$}} }

\newcommand{\bfF}{{\mbox{\boldmath $F$}} }

\newcommand{\bfb}{{\mbox{\boldmath $b$}} }
\newcommand{\bfg}{{\mbox{\boldmath $g$}} }
\newcommand{\bfn}{{\mbox{\boldmath $n$}} }

\newcommand{\half}{\mbox{$\frac{1}{2}$}}

\def\Bbb R{\real}
\def\hat{\widehat}
\def\tilde{\widetilde}
\def\bar{\overline}

\newcommand{\bfchi}{\mbox{\boldmath $\chi$}}

\newcommand{\bfgamma}{\mbox{\boldmath $\gamma$}}

\newcommand{\bfzeta}{\mbox{\boldmath $\zeta$}}

\newcommand{\body}{\mathscr B}
\newcommand{\liquid}{\mathscr L}

\newcommand{\ED}{\end{description}}

\def\tag{\renewcommand{\theequation}}

\newtheorem{theo}{Theorem}
\newtheorem{prop}[theo]{Proposition}
\newtheorem{lem}[theo]{Lemma}
\newtheorem{cor}[theo]{Corollary}
\newtheorem{defi}{Definition}
\newtheorem{rem}{Remark}

\newcommand{\Bd}{\begin{defi}\begin{rm}}
\newcommand{\Ed}{\end{rm}\end{defi}}
\newcommand{\Br}{\begin{rem}\begin{it}}
\newcommand{\Er}{\end{it}\end{rem}}
\newcommand{\Bp}{\begin{prop}\begin{sl}}
\newcommand{\EP}[1]{\end{sl}\label{prop:#1}\end{prop}}
\newcommand{\propref}[1]{{\rm Proposition \ref{prop:#1}}}
\newcommand{\Bt}{\begin{theo}\begin{sl}}
\newcommand{\Et}{\end{sl}\end{theo}}
\newcommand{\Bl}{\begin{lem}\begin{sl}}
\newcommand{\El}{\end{sl}\end{lem}}

\renewcommand{\eqref}[1]{{\rm (\ref{eq:#1})}}
\newcommand{\Bc}{\begin{cor}\begin{sl}}
\newcommand{\Ec}{\end{sl}\end{cor}}
\newcommand{\ET}[1]{\end{sl}\label{theo:#1}\end{theo}}
\newcommand{\EDD}[1]{\end{rm}\label{defi:#1}\end{defi}}
\newcommand{\EL}[1]{\end{sl}\label{lem:#1}\end{lem}}
\newcommand{\theoref}[1]{{\rm Theorem \ref{theo:#1}}}

\newcommand{\ER}[1]{\end{it}\label{rem:#1}\end{rem}}
\newcommand{\EC}[1]{\end{sl}\label{cor:#1}\end{cor}}
\newcommand{\remref}[1]{{\rm Remark \ref{rem:#1}}}
\newcommand{\cororef}[1]{{\rm Corollary \ref{cor:#1}}}
\newcommand{\lemmref}[1]{{\rm Lemma \ref{lem:#1}}}
\newcommand{\essup}[1]{{\rm ess}\,{{\displaystyle \sup_{\hspace*{-5mm}{#1}}}}}

\begin{document}

\title[Stability for a FSI problem]{Stability of equilibria and bifurcations\\
for a fluid-solid interaction problem}
\author{Denis Bonheure
}
\address{
	Département de Mathématique\\
	Université Libre de Bruxelles\\
	Boulevard du Triomphe 155\\
	1050 Brussels - Belgium
}
\email{denis.bonheure@ulb.be}
\date{February 25, 2024}
\author{
Giovanni P. Galdi}\address{
	Dipartment of Mechanical Engineering and Material Sciences\\
	University of Pittsburgh\\
	4200 5th ave, Pittsburgh, PA \\ 15213 USA
	}
 \email{galdi@pitt.edu} \author{Filippo Gazzola}
\address{
	Dipartimento di Matematica\\
	Politecnico di Milano\\
	Piazza Leonardo da Vinci 32\\	
	20133 Milan - Italy}
 \email{filippo.gazzola@polimi.it}

\begin{abstract}
We study certain significant properties of the equilibrium configurations of a rigid body subject to an undamped elastic restoring force, in the stream of a viscous liquid
in an unbounded 3D domain. The motion of the coupled system is driven by a uniform flow at spatial infinity, with constant dimensionless velocity $\lambda$. We show that if $\lambda$ is below a critical value, $\lambda_c$ (say), there is a unique and stable time-independent configuration, where the body is in  equilibrium and the flow is steady. We also prove that, if $\lambda<\lambda_c$, no oscillatory flow may occur. Successively, we investigate  possible loss of uniqueness by providing necessary and sufficient conditions for the occurrence of a steady bifurcation at some $\lambda_s\ge \lambda_c$.


{\bf Primary AMS code.} 76D05, 74F10, 35B32, 76D03
\end{abstract}

\maketitle

\tableofcontents

\renewcommand{\theequation}{\arabic{section}.\arabic{subsection}.\arabic{equation}}

\section{Introduction}\label{chap:intro}

Problems involving the flow of a viscous fluid around solids are the focus of the broad research area of
Fluid-Solid Interactions (FSI). In particular,
also due to their fundamental importance in many practical situations, the oscillations of structures induced by the flow of a viscous liquid occupy a rather  significant position within them. It is thus not surprising that the
problem of flow-induced oscillations  has received all along a plethora of contributions by
the engineering community, from experimental, numerical and theoretical viewpoints; see, e.g., the monographs \cite{Bev,Dyr,PPDL}, the review article \cite{Will} and the references therein. The structure model typically adopted by engineers for this study consists of a rigid body subject to a linear restoring elastic force, while the fluid is modeled by the Navier-Stokes equations \cite{Bev}.\par
Notwithstanding, the problem has not yet received a similar, systematic attention from the mathematical community.
In this regard, in \cite{BBGGP,BGG} we started a rigorous investigation of flow-induced oscillations. There, we dealt with the model problem where a two-dimensional rectangular solid is subject to a unidirectional restoring elastic force,
while immersed in the two-dimensional channel flow of a Navier-Stokes liquid, driven by a time-independent Poiseuille flow.
The main objective concerns the existence and uniqueness of equilibrium configurations of the FSI system, at least for ``small'' flow-rate. Successively, several other works have been dedicated to the investigation of further relevant properties of this model, such as explicit thresholds for uniqueness of the equilibrium configuration \cite{GazP,GazzSp}, non-symmetric configurations \cite{bocchi}, well-posedness of the associated initial-boundary value problem \cite{Patri}, large-time behavior \cite{BoHiPaSpe} and existence of a global attractor \cite{GPP}.\par
Objective of the current paper is to furnish a further contribution to the area of flow-induced oscillations, and consists in the study of the very fundamental properties of the
 stability of equilibrium configurations and possible loss of their uniqueness via steady-state  bifurcations. The model we shall consider --inspired by \cite{Bev}-- is somewhat more general than that in \cite{BGG}, and consists of
 a rigid (finite) body, $\mathscr B$, of {\it arbitrary} shape, subject to a linear undamped restoring force and immersed
in the stream of a Navier--Stokes liquid, $\mathscr L$, that fills the entire three-dimensional space outside $\mathscr B$. The motion of the coupled system is driven by a uniform flow at spatial infinity, characterized by a constant dimensionless velocity $\lambda$ (Reynolds number). The choice of an unbounded (exterior) domain for the flow region is due to the fact that  the interaction of $\mathscr L$ and $\mathscr B$ should not be spoiled by possible ``wall effects".
\par
We are interested in the existence, uniqueness, stability and steady bifurcation of equilibria, where, by
 ``equilibrium'' we mean a state where $\mathscr L$ is in a steady regime and $\mathscr B$ occupies a corresponding fixed region at all times.
\par
Our first goal is to establish the existence of such equilibria, a property that we prove to be valid for {\it all} values of $\lambda$; see \theoref{exi}.  We then provide a variational characterization of their uniqueness by showing the existence of  a critical value $\lambda_1=\lambda_1(\lambda)>0$ such that the
equilibrium is unique if $\lambda-\lambda_1(\lambda)<0$; see \theoref{exi}. As usual, this is merely a sufficient condition for uniqueness which we show to be satisfied (at least) for ``small'' $\lambda$; see \propref{propuniq}.
\par
Successively, we investigate the asymptotic stability (in suitable norms) of the equilibria: our study is complicated by the fact that the region occupied by
the fluid is 3D, unbounded (exterior domain) and, contrary to \cite{BGG}, no Poincar\'e-type inequalities hold. One relevant consequence is that  we are not in a position to furnish a time-decay rate of the perturbations which, very likely, is just algebraic and not exponential as in \cite{GPP}. To set up the stability analysis, we define a second threshold $\lambda_2=\lambda_2(\lambda)\ge0$ such that the stability of the equilibrium is guaranteed if $\lambda-\lambda_2(\lambda)<0$; see \theoref{5.1_01}. However, we are not able to {\it characterize} the equilibria for which the request $\lambda_2(\lambda)>0$ is secured and this condition should be viewed as an assumption which could possibly only hold for certain equilibria. In any case, since $\lambda_2(\lambda)\le\lambda_1(\lambda)$, if $\lambda-\lambda_2(\lambda)<0$ the equilibrium is unique and asymptotically stable. The proof of the latter is carried out by a generalization of the ``invading domains" technique used in \cite{GaSi1}. The main difficulty, in our case, consists in showing that the perturbation to the elongation of the spring eventually tends to 0. Actually, this property is by no means obvious at the outset, since the spring is assumed to be undamped.
\par
Also in view of its importance in the problem of flow-induced oscillations, one may wonder if, in the range $\lambda-\lambda_2(\lambda)<0$,   regimes of oscillatory nature are indeed possible. The answer to this question is given in \theoref{noper}
 where we prove that, in that range of $\lambda$'s, no oscillatory regime can take place. Existence of oscillatory motions can, therefore, take place only for $\lambda>\lambda_2$ as a result of Hopf bifurcation, a question  investigated in the forthcoming article \cite{BoGaGaper}.
\par
The last part the paper is dedicated to steady bifurcation, namely, the existence of possible multiple equilibria for ``sufficiently large'' $\lambda$ (and, certainly, such that $\lambda-\lambda_1(\lambda)\ge 0$). In this regard, we show necessary and sufficient conditions for this phenomenon to occur; see \theoref{2.1.2}.
All these findings are proved by reformulating the equilibrium problem as an operator equation in suitable Banach spaces that allows us to employ known results of abstract bifurcation theory. We emphasize further that, in order to avoid the notorious question of 0 being in the essential spectrum of the linearized operator  \cite{Bab0,FN}, for the functional setting we use {\it homogeneous} (rather than classical) Sobolev spaces, according to the approach introduced in \cite{GaFu}.
\par
The plan of the paper is as follows.
In Section~\ref{chap:prelim} we present the relevant equations and furnish  the mathematical formulation of the problem. In Section \ref{spaces}
we introduce the appropriate functional spaces and collect some of their important properties.
Section \ref{chap:non-osc} contains our main results and is devoted to the study of the equilibria of the FSI problem. Precisely, in Section \ref{sec:exst}
we first prove existence of such equilibria in a class of homogeneous Sobolev spaces, for arbitrary values of the Reynolds number
$\lambda>0$. Successively, we provide the  variational formulation of their uniqueness mentioned above. In the following Section \ref{sec:stab},
we study the asymptotic stability  of  equilibria, whereas the last Section
\ref{steady} is dedicated to the occurrence of steady-state bifurcation.

\section{Formulation of the Problem}\label{chap:prelim}
Let $\mathscr B$ be a rigid body moving in a Navier-Stokes liquid that fills the region $\Omega\subset\mathbb{R}^3$ outside $\mathscr B$ and whose flow becomes uniform at ``large'' distances from $\mathscr B$, characterized by a constant velocity $\bfV\in\mathbb{R}^3$. On $\mathscr B$ an elastic restoring force $\bfF$ acts, applied to its center of mass $G$, while a suitable active torque  prevents it from rotating. Therefore, the motion of $\mathscr B$ is translatory. In this situation, the  governing equations of   motion of the coupled system body-liquid when  referred to a body-fixed frame $\calf\equiv\{G,\bfe_i\}$  are given by \cite[Section 1]{Gah}
\be\ba{cc}\medskip\left.\ba{ll}\medskip
\partial_t\bfv-\nu\Delta\bfv+\nabla p+(\bfv-{\bfgamma})\cdot\nabla\bfv=0\\
\Div\bfv=0\ea\right\}\ \ \mbox{in $\Omega\times(0,\infty)$}\,,\\ \medskip
\bfv(x,t)={\bfgamma}(t)\,, \ \mbox{ $(x,t)\in\partial\Omega\times(0,\infty)$}\,;\ \
\Lim{|x|\to\infty}\bfv(x,t)=\bfV\,,\ t\in(0,\infty)\,,\\
M\dot{\bfgamma}+\rho\Int{\partial\Omega}{} \mathbb T_\nu(\bfv,p)\cdot\bfn=\bfF \ \ \mbox{in $(0,\infty)$}\,.
\ea
\eeq{01}
In \eqref{01}, $\bfv$ and $p$ represent velocity and pressure fields of the liquid,  $\rho$ and $\nu$ its density and kinematic viscosity,  while $M$ and $\bfgamma=\bfgamma(t)$ denote mass of $\mathscr B$  and velocity of $G$, respectively. Here and in the sequel, $\mathbb T_\nu$ denotes the Cauchy stress tensor
$$
\mathbb T_\nu(\bfz,\psi):=2\nu\,\mathbb D(\bfz)-\psi\,\mathbb I\,,\ \ \ \mathbb D(\bfz):=\half\left(\nabla\bfz+(\nabla\bfz)^\top\right)\,,
$$
where $\mathbb I$ is the $3\times3$ identity matrix and $\bfn$ is the unit outer normal at $\partial\Omega$, i.e. directed inside $\body$.\par
We assume that $\bfF$ depends linearly on the displacement $\bfchi(t):=\int\bfgamma(s){\rm d}s=\vec{GO}$, with $O$ fixed point,
namely
\be
\exists\,\ell>0\mbox{ s.t. }\bfF(t)=-\ell\,\bfchi(t),\quad t\ge0.
\eeq{ElFo}
Without loss of generality we take $\bfV=-V\bfe_1$, $V>0$.

\Br
The choice of the linear constitutive equation \eqref{ElFo} is made {\it just} for simplicity of presentation. As will become clear from their proof, our findings (appropriately modified) continue to hold if, more generally, we assume $\bfF=\mathbb A\cdot\bfchi+\bfg(\bfchi)$,   where $\mathbb A$ is a symmetric, positive definite matrix (stiffness matrix), and $\bfg(\bfchi)$ is sufficiently smooth,  with $|\bfg(\bfchi)|=o(|\bfchi|)$ as $|\bfchi|\to 0$.
\ER{0.0}

Writing $\bfv=\bfu-V\bfe_1$, we are led to
\be\ba{cc}\medskip\left.\ba{ll}\medskip
\partial_t\bfu-\nu\Delta\bfu+\nabla p+(\bfu-{\dot\bfchi}(t))\cdot\nabla\bfu-V\bfe_1\cdot\nabla \bfu=0\\
\Div\bfu=0\ea\right\}\ \ \mbox{in $\Omega\times(0,\infty)$}\,,\\ \medskip
\bfu(x,t)={\dot\bfchi}(t)\!+\!V\bfe_1,\ (x,t)\in\partial\Omega\times(0,\infty);\
\Lim{|x|\to\infty}\bfu(x,t)=0,\ t\in(0,\infty),\\
M\ddot{\bfchi}+\ell\bfchi+\rho\Int{\partial\Omega}{} \mathbb T_\nu(\bfu,p)\cdot\bfn=0 \ \ \mbox{in $(0,\infty)$}\,.
\ea
\eeq{01bis}

Scaling  velocity with $V$, length with $L:={\rm diam}\,\mathscr B$,  time with $L^2/\nu$,
and setting $\bfu:=\bfv+\bfe_1$, we may rewrite \eqref{01bis} in the following form
\be\ba{cc}\medskip\left.\ba{ll}\medskip
\partial_t\bfu-\Delta\bfu+\nabla p=\lambda\,[\partial_1\bfu+(\dot{\bfchi}-\bfu)\cdot\nabla\bfu]\\
\Div\bfu=0\ea\right\}\ \ \mbox{in $\Omega\times(0,\infty)$}\,,\\ \medskip
\bfu(x,t)={\dot{\bfchi}}(t)+\bfe_1\,, \ \mbox{ $(x,t)\in\partial\Omega\times(0,\infty)$}\,;\ \
\Lim{|x|\to\infty}\bfu(x,t)=\0\,,\ t\in(0,\infty)\,,\\
\ddot{\bfchi}+\omega_{\sf n}^2\bfchi+\varpi\Int{\partial\Omega}{} \mathbb T_1(\bfu,p)\cdot\bfn=\0 \ \ \mbox{in $(0,\infty)$}\,,
\ea
\eeq{02}
with
$$
\omega_{\sf n}^2:=\frac{L^4\ell }{M\nu^2}\,,\    \ \varpi:=\frac{\rho L^3}{M}\,,\ \ \lambda:=\frac{VL}{\nu}\,.
$$%
All the involved quantities are now non-dimensional; in the sequel, we just write $\mathbb T$ instead of $\mathbb T_1$.

Let ${\sf s}_0=(\bfu_0,p_0,\bfchi_0)$ be a steady-state solution to \eqref{02} corresponding to a given $\lambda$, namely,
\be\ba{cc}\medskip\left.\ba{ll}\medskip
-\Delta\bfu_0+\nabla p_0=\lambda\,(\partial_1\bfu_0-\bfu_0\cdot\nabla\bfu_0)\\
\Div\bfu_0=0\ea\right\}\ \ \mbox{in $\Omega$}\,,\\ \medskip
\bfu_0(x)=\bfe_1\,, \ \mbox{ $x\in\partial\Omega$}\,;\ \
\Lim{|x|\to\infty}\bfu_0(x)=\0\,,\\
\omega_{\sf n}^2\bfchi_0+\varpi\Int{\partial\Omega}{} \mathbb T(\bfu_0,p_0)\cdot\bfn=\0\,.
\ea
\eeq{03}
From the physical viewpoint, $\bfchi_0$ represents the (non-dimensional and rescaled) elongation of the spring necessary to keep $\calb$ in place.
In Sections \ref{sec:exst} and \ref{sec:stab}, for any $\lambda>0$, we show that \eqref{03} has at least one solution
$${\sf s}_0(\lambda):=(\bfu_0(\lambda),p_0(\lambda),\bfchi_0(\lambda))$$
that is unique and stable provided $\lambda$ remains  below a definite value $\lambda_c$ that we characterize in \theoref{5.1_01}. The uniqueness threshold does only depend on $\Omega$, see \propref{propuniq}, whereas the stability threshold $\lambda_c$ depends on the solution itself and basically on its decay at spatial infinity, see \remref{max}.
 Moreover, we prove that as long as $\lambda<\lambda_c$, no oscillatory regime can branch out of ${\sf s}_0(\lambda)$; see Section~\ref{sec:nobif}.
Therefore, non-uniqueness of, as well as bifurcation from ${\sf s}_0(\lambda)$ may occur only at some $\lambda \ge \lambda_c$. In this regard,
we next investigate the occurrence of steady-state bifurcation  at some $\lambda=\lambda_s\ge\lambda_c$. More precisely, we furnish necessary and sufficient conditions  for the existence of a bifurcation point $\lambda_s$ and a family of solutions ${\sf s}(\lambda)$ to \eqref{03}, with ${\sf s}(\lambda) \not\equiv {\sf s}_0(\lambda)$,  $\lambda\in U(\lambda_s)$, such that ${\sf s}(\lambda)\to{\sf s}_0(\lambda_s)$, as $\lambda\to\lambda_s$. This is accomplished by formulating the problem in a functional setting that is suitable to employ the abstract results of bifurcation theory; see Section~\ref{steady}.

\setcounter{equation}{0}
\section{Functional framework}\label{spaces}

\subsection{Notations and Relevant Functional Spaces}
Let $\Omega_0\subset\mathbb{R}^3$ be the closure of a bounded domain of class $C^2$, representing the region occupied by $\mathscr B$.
Let $\Omega=\mathbb{R}^3\setminus\Omega_0$ be the  unbounded exterior domain containing the fluid $\liquid$. With the origin of coordinates in
the interior of $\Omega_0$, we set
$$
B_R:=\{x\in\real^3:\,|x|<R\},\ \Omega_R:=\Omega\cap B_R,\
\Omega^R:=\Omega\backslash\bar{\Omega_R}\quad\forall R>R_*:={\rm diam}\, \Omega_0.$$

As customary, for a domain
$A\subset\real^3$,  $L^q=L^q(A)$ denotes the Lebesgue space with norm $\|\cdot\|_{q,A}$, and  $W^{m,2}=W^{m,2}(A)$, $m\in\nat$,  the Sobolev space with norm $\|\cdot\|_{m,2,A}$. By $(\,\ ,\,\ )_A$ we indicate the $L^2(A)$-scalar product.  Furthermore, $D^{m,q}=D^{m,q}(A)$ is the homogeneous Sobolev space with semi-norm $\sum_{|l|=m}\|D^lu\|_{q,A}$, whereas $D_0^{1,2}=D_0^{1,2}(A)$ is the completion of $C_0^\infty(A)$ in the norm $\|\nabla(\cdot)\|_{2,A}$. The dual space of $D_0^{1,2}(A)$ will be indicated by $D_0^{-1,2}(A)$. In all the above notation we shall typically omit the subscript ``$A$'', unless confusion arises.\par
If $M$ is a map between two Banach spaces $X$ and $Y$, we denote by ${\sf D}[M]\subseteq X$ and ${\sf R}[M]\subseteq Y$ its domain and range, respectively,  and by ${\sf N}[M]:=\{x\in X: M(x)=0\}$ its null space.\par
We shall now introduce  certain function classes characterized by the property that their elements are solenoidal. Their most important  properties will be collected later on in Section \ref{sub:24}. Let
$$\ba{ll}\medskip
\calk=\mathcal K(\real^3):=\big\{\bfphi\in C_0^\infty({\real^3}):
\exists\, \hat{\bfphi}\in\real^3 \text{ s.t. }\bfphi(x)\equiv\hat{\bfphi} \text{ in a neighborhood of }\Omega_0\big\}
\\ \medskip
\calc=\mathcal C(\real^3):=\{\bfphi\in\calk(\real^3):\ \Div\bfphi=0\ \mbox{in $\real^3$}\}\,,\\
\calc_0=\mathcal C_0(\Omega):=\{\bfphi\in\calc(\real^3): \hat{\bfphi}=\0\}\,.
\ea
$$

In $\calk$ we consider the scalar product
\be
\langle \bfphi,\bfpsi\rangle:=
\varpi^{-1}\,\hat{\bfphi}\cdot\hat{\bfpsi}+(\bfphi,\bfpsi)_\Omega\,,\quad\forall\bfphi,\bfpsi\in\calk\,,
\eeq{0.0}
and we introduce the spaces
\arraycolsep=1pt
\be\ba{rl}
\mathcal L^2=\call^2(\real^3):= &\left\{
\text{completion of }\calk(\real^3)
\text{ in the norm induced by }\eqref{0.0}\right\}\\
\mathcal H=\calh(\real^3):= &\left\{\text{completion of }\calc(\real^3) \text{ in the norm induced by \eqref{0.0}}\right\} \\
\calg=\calg(\real^3):=
&\big\{\bfh\in \call^2(\real^3) : \, \exists\, p\in D^{1,2}(\Omega) \text{ s.t. }
\bfh=\nabla p \text{ in } \Omega,
\\
& \text{and } \bfh=-\varpi\int_{\partial\Omega}p\,\bfn \text{ in } \Omega_0\big\}.
\ea
\eeq{spazi}
\arraycolsep=5pt
\par
We next define
$$\ba{ll}\medskip
\cald^{1,2}=\cald^{1,2}(\real^3):=\,\big\{\mbox{completion of $\calc(\real^3)$ in the norm $\|\mathbb D(\cdot)\|_2$}\big\}\,,\\
\cald_0^{1,2}=\cald_0^{1,2}(\Omega):=\,\big\{\mbox{completion of $\calc_0(\Omega)$ in the norm $\|\mathbb D(\cdot)\|_2$}\big\}\,,
\ea
$$
\be
Z^{2,2}:=W^{2,2}(\Omega)\cap \cald^{1,2}(\real^3)\,.
\eeq{zetaa}
Obviously, $\cald_0^{1,2}(\Omega)\subset\cald^{1,2}(\real^3)$.\par
Along with the spaces $\call^2,\calh$, and  $\cald^{1,2}$ defined above, we shall need also their ``restrictions'' to the ball $B_R$. Precisely, we set
$$\ba{ll}\medskip
\call^2(B_R):= \{\bfphi\in L^2(B_R): \, \bfphi|_{\Omega_0}=\hat{\bfphi}\,\ \mbox{for some $\hat{\bfphi}\in\real^3$}\}\\ \medskip
\calh(B_R):=\{\bfphi\in \call^2(B_R):\, \Div\bfphi=0\,,\ \bfphi\cdot\bfn|_{\partial B_R}=0\}\\
\cald^{1,2}(B_R):= \{\bfphi\in W^{1,2}(B_R)\cap \call^2(B_R) : \, \Div\bfphi=0\,,\ \bfphi|_{\partial B_R}=\0 \}\,.
\ea
$$
Then $\calh(B_R)$ and $\cald^{1,2}(B_R)$ are Hilbert spaces with scalar products
\[
\begin{split}
\varpi^{-1}\,\hat{\bfphi}_1\cdot\hat{\bfphi}_2+(\bfphi_1,\bfphi_2)_{\Omega_R}\,,\ \bfphi_i\in\calh(\Omega_R)\,;\ \\
(\mathbb D(\bfpsi_1),\mathbb D(\bfpsi_2))\,,\ \bfpsi_i\in \cald^{1,2}(\Omega_R)\,,\ \ i=1,2.
 \end{split}
\]

Let $\cald_0^{-1,2}(\Omega)$ be the dual space of $\cald_0^{1,2}(\Omega)$, endowed with the norm
$$
|\bff|_{-1,2}=\Sup{\mbox{\footnotesize $\ba{c}\bfphi\in \calc_0(\Omega)\\ \|\nabla\bfphi\|_2=1\ea $}}\left|\left(\bff,\bfphi\right)\right|\,,
$$
and set
\be
\caly:=\cald_0^{-1,2}(\Omega)\cap \calh(\real^3)\,,\ \ Y:=\cald_0^{-1,2}(\Omega)\cap \call^2(\real^3)\,,
\eeq{ips}
with associated norms
$$
\|\bfg\|_{\caly}=\|\bfg\|_{Y}:=\|\bfg\|_{2}+|\bfg|_{-1,2}+|\hat{\bfg}|\,.
$$
We then define
\be
\begin{split}
X=X(\Omega):=\{\bfu\in\cald^{1,2}_0(\Omega):\, \partial_1\bfu\in \cald_0^{-1,2}(\Omega)\},\\ X^2=X^2(\Omega):=\left\{\bfu\in X(\Omega): D^2\bfu\in L^2(\Omega)\right\}\,.
\end{split}
\eeq{ics}
It is known \cite[Proposition 65]{GaCe} that $X$ and $X^2$ are (reflexive, separable) Banach spaces when equipped with the norms
$$
\|\bfu\|_X:=\|\nabla\bfu\|_2+|\partial_1\bfu|_{-1,2}\,,\qquad
\|\bfu\|_{X^2}:=\|\bfu\|_{X}+\|D^2\bfu\|_{2}\,.
$$
In fact, as shown later on in \lemmref{1.1}, the norms $\|\nabla(\cdot)\|_2$ and $\|\mathbb D(\cdot)\|_2$ are equivalent in $\cald^{1,2}_0$.\par
In the sequel, we also need some spaces of {\it time-periodic functions}. A function $\bfw:\Omega\times \real\mapsto \real^3$ is
{\it $2\pi$-periodic}, if for a.e.\ $t\in \real$, $\bfw(\cdot,t+~2\pi)=\bfw(\cdot,t)$,
and we use the standard notation
\be
{\bar \bfw(\cdot)}:=\frac{1}{2\pi}\int_{0}^{2\pi}\bfw(\cdot,t){\rm d}t\,,
\eeq{media}
whenever the integral is meaningful.
Let $B$ be a function space  with seminorm $\|\cdot\|_B$. By $L^2(0,2\pi;B)$ we denote the class of functions
$u:(0,2\pi)\rightarrow B$ such that
$$
\|u\|_{L^2(B)}:= \left(\Int{0}{2\pi}\|u(t)\|_B^2 {\rm d}t\right)^{\frac 12}<\infty
$$
Likewise, we put
$$
W^{1,2}(0,2\pi;B)=\Big\{u\in L^{2}(0,2\pi;B): \partial_t u\in L^{2}(0,2\pi;B)\Big\}\,.
$$
For simplicity, we write $L^2(B)$ for $L^2(0,2\pi;B)$, etc. Moreover, we define the Banach spaces
$$\ba{rl}\medskip
L^{q}_\sharp&\!\!\!\!:=\{\bfxi\in L^{q}(0,2\pi), \ \mbox{$\bfxi$ is $2\pi$-periodic with $\bar{\bfxi}=\0$}\}\,,\ \ q\in[1,\infty]\,,\\\medskip
W^{k}_\sharp&\!\!\!\!:=\{\bfxi\in L^{2}_\sharp(0,2\pi), \ d^l{\bfxi}/dt^l\in L^2(0,2\pi)\,,\ \ l=1,\ldots,k\}\,,\\\medskip
\mathcal L_\sharp^{2}&\!\!\!\!:=\{\bfw\in L^{2}(L^2(\Omega)); \ \mbox{$\bfw$ is $2\pi$-periodic, with $\bar{\bfw}=\0$}\}\,,
\\
\mathcal W_\sharp^{2}&\!\!\!\!:=\{\bfw\in W^{1,2}(L^2(\Omega))\cap L^2(W^{2,2}(\Omega)); \ \mbox{$\bfw$ is $2\pi$-periodic, with $\bar{\bfw}=\0$}\}\,,
\ea
$$
with associated norms
$$\ba{ll}\medskip
\|\bfxi\|_{L^q_\sharp}:=\|\bfxi\|_{L^q(0,2\pi)}\,,\ \ \ \|\bfxi\|_{W^k_\sharp}:=\|\bfxi\|_{W^{k,2}(0,2\pi)}\,,\\
\|\bfw\|_{\mathcal L^{2}_\sharp}:=\|\bfw\|_{L^2(L^2(\Omega))}\,,\ \ \
\|\bfw\|_{\mathcal W_\sharp^{2}}:=\|\bfw\|_{W^{1,2}(L^2(\Omega))}+\|\bfw\|_{L^2(W^{2,2}(\Omega))}\,.
\ea
$$
We also introduce the Banach spaces
$$\ba{rl}\medskip
\textsf{{W}$_\sharp^{2}$}&\!\!\!\!:=\left\{\bfw\in L^2(Z^{2,2})\cap W^{1,2}(\calh):\, \mbox{$\bfw$ is $2\pi$-periodic,\,with\, $\bar{\bfw}|_{\Omega_0}=\bar{\hat{\bfw}}=\0$}\right\}\,,
\\
\textsf{{L}$_\sharp^{2}$}&\!\!\!\!:=\left\{\bfw\in L^2(\calh):\, \mbox{$\bfw$ is $2\pi$-periodic,\,with\, $\bar{\bfw}|_{\Omega_0}=\bar{\hat{\bfw}}=\0$}\right\}
\ea
$$
with corresponding norms
$$
\|\bfw\|_{\mbox{\scriptsize $\W$}}:=\|\partial_t\bfw\|_{L^2(\Omega)}+\|\bfw\|_{L^2(W^{2,2}(\Omega))}+\|\hat{\bfw}\|_{W_\sharp^1}\,,\qquad
\|\bfw\|_{\mbox{\scriptsize $\M$}}:=\|\bfw\|_{L^2(L^2(\Omega))}+\|\hat{\bfw}\|_{L^2_\sharp}\,.
$$

Finally, we set
$$
{\sf P}^{1,2}_\sharp:=\left\{{\sf p}\in L^2(D^{1,2})\ \mbox{with \ $\bar{\sf p}=0$}\right\}\,,
$$
with associated norm
$$
\|{\sf p}\|_{{\sf P}^{1,2}_\sharp}:=\|{\sf p}\|_{ L^2(D^{1,2})}\,.
$$

\subsection{Basic Properties of the Relevant Functional Spaces}\label{sub:24}

The $\call^2$-spaces have two main properties.

\Bl The following characterizations hold
$$
\call^2(\real^3)=\{\bfu\in L^2(\real^3): \ \bfu=\hat{\bfu}\ \mbox{in}\ \Omega_0, \text{ for some }\hat{\bfu}\in \real^3\},\quad
\calh(\real^3)=\{\bfu\in \call^2(\real^3): \ \Div\bfu=0\,\}\,.
$$
\EL{aria}
{\it Proof.} See  \cite[Theorem 3.1 and Lemma 3.2]{ALS}. \hfill$\square$\par\smallskip

We also have
\Bl With the scalar product \eqref{0.0}, the following orthogonal decomposition holds
$$
\call^2(\real^3)=\calh(\real^3)\oplus\calg(\real^3)\,.
$$
\EL{0.1}
{\it Proof.} A proof can be deduced from  \cite[Theorem 3.2]{ALS}. However, for completeness and since this result plays a major role in our analysis, we reproduce it here. Let $\bfu\in \calh$ and $\bfh\in\calg$. Then,
$$
\langle\bfu,\bfh\rangle=\int_{\Omega}\bfu\cdot\nabla p-\int_{\partial\Omega}p\,\hat{\bfu}\cdot\bfn\,.
$$
Therefore, integrating by parts and using $\Div\bfu=0$ we deduce
$$
\langle\bfu,\bfh\rangle=-\int_{\Omega}p\,\Div\bfu+\int_{\partial\Omega}p\,\hat{\bfu}\cdot\bfn-\int_{\partial\Omega}p\,\hat{\bfu}\cdot\bfn=0
$$
which proves $\calh^{\perp}\supset \calg$. Conversely, assume  $\bfv\in\calh^{\perp}$, i.e.
\be
\varpi^{-1}\hat{\bfv}\cdot\hat{\bfu}+\int_\Omega\bfv\cdot\bfu=0\,,\ \ \mbox{for all $\bfu\in\calh$.}
\eeq{dip}
Since $\calc_0\subset \calh$, by picking $\bfu\in \calc_0$ from the preceding  we find, in particular,
$$
\int_\Omega\bfv\cdot\bfu=0\,,\ \ \mbox{for all $\bfu\in\calc_0$,}
$$
so that, by well-known results on the Helmholtz decomposition \cite[Lemma III.1.1]{Gab}, we infer $\bfv=\nabla p$ with $p\in D^{1,2}(\Omega)$. Replacing the latter into \eqref{dip} and integrating by parts, we get
$$
\left(\varpi^{-1}\hat{\bfv}+\int_{\partial\Omega}p\,\bfn\right)\cdot\hat{\bfu}=0\ \ \mbox{for all $\hat{\bfu}\in\real^3$}\,,
$$
from which we conclude that $\bfv\in\calg$, that is, $\calh^\perp\subset \calg$. The proof of the lemma is completed.\hfill$\square$\par\smallskip

Concerning the properties of $\cald^{1,2}$-spaces, we state

\Bl  Let $\tilde{\cald^{1,2}}$ denote either $\cald^{1,2}$ or $\cald_0^{1,2}$. Then, $\tilde{\cald^{1,2}}$ is a separable Hilbert space when equipped with the scalar product
$$
(\mathbb D(\bfu_1),\mathbb D(\bfu_2))\,,\ \ \bfu_i\in \tilde{\cald^{1,2}}\,, \, \ i=1,2\,.
$$
Moreover, we have the characterization:
\be
\tilde{\cald^{1,2}}=\big\{\bfu\in L^6(\real^3)\cap D^{1,2}(\real^3)\,;\ \Div\bfu=0\,;\,
\bfu=\hat{\bfu} \ \mbox{in $\Omega_0$} \big\}\,,
\eeq{1.7}
with some $\hat{\bfu}\in\real^3$ if $\tilde{\cald^{1,2}}\equiv\cald^{1,2}$, and $\hat{\bfu}=\0$ if $\tilde{\cald^{1,2}}\equiv\cald_0^{1,2}$.
Also, for each $\bfu\in\tilde{\cald^{1,2}}$, it holds
\be
\|\nabla\bfu\|_2=\sqrt{2}\|\mathbb D(\bfu)\|_2\,,
\eeq{1.8}
and
\be
\|\bfu\|_6\le \kappa_0\,\|\mathbb D(\bfu)\|_2\,,
\eeq{1.9}
for some $\kappa_0>0$. Finally,
there is another positive constant $\kappa_1$ such that
\be
|\hat{\bfu}|\le \kappa_1\,\|\mathbb D(\bfu)\|_2\,.
\eeq{1.10}
\EL{1.1}
{\it Proof.} See \cite[Lemmas 9--11]{Gah}.\hfill$\square$\par\smallskip

\Br The space $\cald^{1,2}(B_R)$ can be viewed as a subspace of $\cald^{1,2}(\real^3)$, by extending to 0 in $\real^3\backslash B_R$ its generic element. Therefore, all the properties mentioned in \lemmref{1.1} continue to hold for $\cald^{1,2}(B_R)$.
\ER{2.1}

The $X$-spaces also have a number of relevant properties that we collect in the next statements.

\Bl The following continuous embedding properties hold
\be
X^2(\Omega)\subset W^{2,2}(\Omega_R)\ \ \mbox{for all $R>R_*$}\,,\quad
X^2(\Omega)\subset  L^\infty(\Omega)\cap D^{1,q}(\Omega)\ \ \mbox{for all $q\in[2,6]$}.
\eeq{X2e}\EL{Xemb}
{\it Proof.} By \lemmref{1.1}, the first property is obvious. From \cite[Theorem II.6.1(i)]{Gab} it follows that $X^2(\Omega)\subset D^{1,6}(\Omega)$ which, in turn, by \cite[Theorem II.9.1]{Gab} and simple interpolation allows us to deduce
also the second stated property.\hfill$\square$\par\smallskip
We conclude this section with the following embedding result, whose proof is given in  \cite[Lemma 2]{GaNe}.
\Bl The space $X(\Omega)$ is continuously embedded in $L^4(\Omega)$.
\EL{3.1}

\setcounter{equation}{0}
\section{Main results on the equilibrium configurations}\label{chap:non-osc}

\subsection{Existence and Uniqueness}\label{sec:exst}
We begin with a general existence result  in a suitable function class, followed by a  corresponding uniqueness result. Both findings are, in fact, obtained as a corollary to classical results regarding steady-state Navier-Stokes problems in exterior domains.
Precisely, we have the following theorem.
\Bt For any $\lambda>0$, problem \eqref{03} has at least one solution
$${\sf s}_0(\lambda):=(\bfu_0(\lambda),p_0(\lambda),\bfchi_0(\lambda))$$ such that
\be {\sf s}_0(\lambda)\in [L^q(\Omega)\cap D^{1,r}(\Omega)\cap D^{2,s}(\Omega)]\times[L^{\sigma}(\Omega)\cap D^{1,s}(\Omega)]\times\real^3,
\eeq{sfaco}
for all $q\in (2,\infty]$, $r\in(\frac43,\infty]$, $s\in (\frac32,\infty]$, $\sigma\in (1,\infty)$.  Moreover, the quantity
\be \sup_{{\mbox{\footnotesize $\bfu\in\cald^{1,2}_0(\Omega)$}}}\frac{(\bfu\cdot\nabla\bfu,\bfu_0)}{\|\nabla\bfu\|_2^2}=:\frac1{\lambda_1}
\eeq{sfaco1}
is finite, strictly positive, and achieved and, if $\lambda<\lambda_1$,  the solution ${\sf s}_0(\lambda)$ is unique.
\medbreak

\ET{exi}
{\it Proof.} From \cite[Theorem X.6.4]{Gab} we know that for any $\lambda>0$  problem \eqref{03}$_{1-4}$ has  one corresponding  solution $(\bfu_0,p_0)$ in the class \eqref{sfaco}. We then set
\be
\bfchi_0:= -\frac{\varpi}{\omega_{\sf n}^2}\int_{\partial\Omega}\mathbb T(\bfu_0,p_0)\cdot\bfn\,,
\eeq{2.3_00}
which is well defined by standard trace theorems. This completes the proof of the existence.

We now turn to the uniqueness part. The existence and achievement of $1/\lambda_1$ follows from the summability properties of $\bfu_0$ given in \eqref{sfaco} and standard arguments about maxima of quadratic forms in exterior domains \cite{Gam}. In order to prove that $\lambda_1>0$, take any
$\bfw\in C^\infty_0(\real^3)$ with $\supp\bfw\cap\Omega_0=\emptyset$. Then let $\bfu=\curl\bfw$ so that $\bfu\in\cald^{1,2}_0(\Omega)$.
By translating rigidly $\bfu$ and moving its support towards infinity, we see that
$(\bfu\cdot\nabla\bfu,\bfu_0)\to0$ due to the decay properties of $\bfu_0$. Therefore, by its characterization in \eqref{sfaco1},
$\lambda_1^{-1}\ge0$; in fact, $\lambda_1^{-1}>0$ since the supremum in \eqref{sfaco1} is achieved.\par
Finally, let $(\bfu+\bfu_0,p+p_0,\bfchi+\bfchi_0)$ be another solution to \eqref{03} in the class \eqref{sfaco} corresponding to the same $\lambda$. Then
$(\bfu,p,\bfchi)$ satisfies the following equations
\be\ba{cc}\medskip\left.\ba{ll}\medskip
-\Delta\bfu+\nabla p=\lambda\,[\partial_1\bfu-(\bfu_0+\bfu)\cdot\nabla\bfu-\bfu\cdot\nabla\bfu_0]\\
\Div\bfu=0\ea\right\}\ \ \mbox{in $\Omega$}\,,\\
\medskip
\bfu(x)=\0\,, \ \mbox{ $x\in\partial\Omega$}\,,\qquad \Lim{|x|\to\infty}\bfu(x)=\0 \,,\\
\omega_{\sf n}^2\bfchi+\varpi\Int{\partial\Omega}{} \mathbb T(\bfu,p)\cdot\bfn=\0\,.
\ea
\eeq{03_111}
Dot-multiplying \eqref{03_111}$_1$ by $\bfu$, integrating by parts over $\Omega$ and using \eqref{sfaco}, \eqref{03_111}$_{2,3}$ and \eqref{sfaco1} we find
$$
\|\nabla\bfu\|_2^2=\lambda(\bfu\cdot\nabla\bfu,\bfu_0)\le \frac{\lambda}{\lambda_1}\|\nabla\bfu\|_2^2\,,
$$
from which it follows that $\lambda\ge\lambda_1$ or $\|\nabla\bfu\|_2=0$. The latter implies $\bfu=\0$, $\bfchi=\0$ (and therefore $p=0$) by \lemmref{1.1}.
\hfill$\square$\par\smallskip

\Br
The value $\lambda_1$ defined in \eqref{sfaco1} plays the role of a weighted Poincar\'e constant. Indeed, it may be equivalently
characterized by
$$
\lambda_1= \min_{{\mbox{\footnotesize $\bfu\in\cald^{1,2}_0(\Omega)$}}}\frac{\|\nabla\bfu\|_2^2}{(\bfu\cdot\nabla\bfu,\bfu_0)}
=\min_{{\mbox{\footnotesize $\bfu\in\cald^{1,2}_0(\Omega)$}}}\frac{\|\nabla\bfu\|_2^2}{(\bfu\cdot\nabla(-\bfu_0),\bfu)}
$$
with the weight $\nabla(-\bfu_0)$ vanishing at infinity and bringing enough compactness to ensure that the minimum is achieved.
\Er

Since $\lambda_1$ depends on $\bfu_0$ which in turn depends on $\lambda$, it is natural to wonder whether the condition $\lambda<\lambda_1$
(ensuring uniqueness) can be reached. The next statement shows that this is the case.

\Bp There exists $\gamma=\gamma(\Omega)>0$ such that, if $\lambda<\gamma$, then problem \eqref{03} admits a unique solution.
\medbreak

\EP{propuniq}

{\it Proof.} As already noticed in the existence proof, the fluid equations \eqref{03}$_{1-4}$ decouple from the one in \eqref{03}$_5$, representing the balance of forces on $\mathscr B$. Therefore, uniqueness for the whole problem  \eqref{03} is reduced to establish the same property just for the Navier-Stokes problem \eqref{03}$_{1-4}$. However, the latter is well known \cite[Theorem X.7.3]{Gab} and is achieved exactly under the condition stated in the proposition.
\par\hfill$
\square$\par

From \propref{propuniq} and \theoref{exi}, we infer, in particular, that uniqueness is ensured for ``small" $\lambda>0$ and may fail only at some $\lambda$ such that
\be
\frac1\lambda\le\max_{\bfu\in\cald^{1,2}_0(\Omega)}\frac{(\bfu\cdot\nabla\bfu,\bfu_0(\lambda))}{\|\nabla\bfu\|_2^2}.
\eeq{inequall}

This reveals that either we have uniqueness for all $\lambda>0$ or there exists $\tilde\lambda>0$ such that $\lambda_1(\tilde\lambda)=\tilde\lambda$.
This $\tilde\lambda\in(0,\infty]$ can be defined as
$$\tilde\lambda=\sup\{\lambda>0: \nu<\lambda_1(\nu),\forall \nu\in (0,\lambda)\}.$$
If $\tilde\lambda<\infty$, then there exists a non trivial solution $\bfu$ of the linear equation
\be\ba{cc}\medskip\left.\ba{ll}\medskip
\Delta\bfu-\nabla p = \tilde\lambda\,(\bfu_0(\tilde\lambda)\cdot\nabla\bfu+\bfu\cdot\nabla\bfu_0(\tilde\lambda))\\
\Div\bfu=0\ea\right\}\ \ \mbox{in $\Omega$}\,,\\
\medskip
\bfu(x)=\0\,, \ \mbox{ $x\in\partial\Omega$}\,,  \Lim{|x|\to\infty}\bfu(x)=\0\,.
\ea
\eeq{maybe}
This condition is, in general, only necessary to get multiple equilibria (i.e.\ for $\lambda>\tilde\lambda$), as discussed in detail later on in the bifurcation context; see \theoref{2.1} and \cororef{6.1}.

From a physical point of view, one expects that $\bfu_0=\bfu_0(\lambda)$ becomes ``larger'' as $\lambda$ grows, although a precise
definition of ``larger'' appears out of reach. From a mathematical point of view, this could be translated into the fact that some norms of
$\bfu_0(\lambda)$ are expected to grow with $\lambda$. If this were true, then equality in \eqref{inequall} would hold for a unique value
$\overline{\lambda}>0$ and this would imply that
$$
\lambda_1(\lambda)>\lambda\mbox{ if }\lambda<\overline{\lambda}\, ,\quad\lambda_1(\lambda)<\lambda\mbox{ if }\lambda>\overline{\lambda}\, .
$$
Clearly, this {\it would not} allow us to conclude that uniqueness for \eqref{03} is ensured if and only if $\lambda<\overline{\lambda}$.

\setcounter{equation}{0}
\subsection{Asymptotic stability}
\label{sec:stab}

\subsubsection{A sufficient condition for stability}

Our next task is to find sufficient conditions for the stability of solutions determined in \theoref{exi}, in a suitable class of ``perturbations''. In this regard, let ${\sf s}_0(\lambda)=(\bfu_0(\lambda),p_0(\lambda),\bfchi_0(\lambda))$ be a steady-state solution
given in \eqref{sfaco} and let $(\bfu,p,\bfchi)$ be a corresponding time-dependent perturbation. By \eqref{02}, we then have that $(\bfu,p,\bfchi)$ satisfies the following set of equations
\be\ba{c}\medskip\left.\ba{r}\medskip
\partial_t\bfu-\Delta\bfu+\nabla {p}
 =\lambda\,[\partial_1\bfu-\bfu_0\cdot\nabla\bfu+(\dot{\bfchi}-\bfu)\cdot\nabla\bfu_0 -(\dot{\bfchi}-\bfu)\cdot\nabla\bfu]\\
\Div\bfu=0\ea\right\}\ \ \mbox{in $\Omega\times(0,\infty)$}\,,\\ \medskip
\bfu(x,t)={\dot{\bfchi}}(t)\,, \ \mbox{ $(x,t)\in\partial\Omega\times(0,\infty)$}\,,\\ \medskip
\ddot{\bfchi}+\omega_{\sf n}^2\bfchi+\varpi\Int{\partial\Omega}{} \mathbb T(\bfu,{p})\cdot\bfn=\0\,, \ \ \mbox{in $(0,\infty)$}\\
\bfu(x,0)=\bfu^0\,,\ x\in\Omega\,,\ \ \bfchi(0)=\bfchi^{0}\,,\ \ \dot{\bfchi}(0)=\bfchi^1
\,.
\ea
\eeq{04_111}

Given $\lambda>0$ and some
${\sf s}_0(\lambda)$ solving problem \eqref{03} (see \theoref{exi}) we define
\be
\frac1{\lambda_2}=\frac1{\lambda_2(\lambda)}:=\sup_{\mbox{\footnotesize $\ \bfu\in\cald^{1,2}(\real^3)$}}\frac{((\bfu-\hat{\bfu})\cdot\nabla\bfu,\bfu_0(\lambda))}{\|\nabla\bfu\|_2^2}
\eeq{sfaco2}
From \theoref{exi} we know that
$$\frac1{\lambda_1}=\frac{(\tilde{\bfu}\cdot\nabla\tilde{\bfu},\bfu_0)}{\|\nabla\tilde{\bfu}\|_2^2},$$
for some $\tilde{\bfu}\in \cald_0^{1,2}(\Omega)$. Moreover, as $\cald_0^{1,2}(\Omega)\subset \cald^{1,2}(\real^3)$
(see Section~\ref{spaces}), and since $\tilde{\bfu}$ vanishes on $\partial\Omega$, we infer that
$$\frac1{\lambda_1}=\frac{(\tilde{\bfu}\cdot\nabla\tilde{\bfu},\bfu_0)}{\|\nabla\tilde{\bfu}\|_2^2}\le \sup_{\mbox{\footnotesize $\ \bfu\in\cald^{1,2}(\real^3)$}}\frac{((\bfu-\hat{\bfu})\cdot\nabla\bfu,\bfu_0)}{\|\nabla\bfu\|_2^2}=\frac1{\lambda_2}.$$
Therefore, we infer
\be
0\le\lambda_2\le \lambda_1\, .
\eeq{lambda2ge0}

We can now state a stability result for \eqref{04_111}.

\Bt Let $\lambda>0$ and let ${\sf s}_0(\lambda)$ be a solution of problem \eqref{03}. Suppose $\bfu_0(\lambda)$ is such that $\lambda_2^{-1}(\lambda)<\infty$ and that $\lambda<\lambda_2$.
Then, there exists  $\varepsilon=\varepsilon(\Omega,\lambda,\omega_{\sf n},\varpi)>0$ such that, if
\be
\|\bfu^0\|_{1,2}+|\bfchi^0|+|\bfchi^1|\le\varepsilon,
\eeq{smalleps}
then problem \eqref{04_111} has one and only one solution such that
\be\left.\ba{ll}\medskip
\bfu\in C([0,T]; \cald^{1,2}(\real^3))\cap L^2(0,T; W^{2,2}(\Omega))\cap W^{1,2}(0,T; \call^2(\real^3))\,,\\ p\in L^2(0,T;D^{1,2}(\Omega))\,,\ \ \bfchi\in W^{2,2}(0,T)\,,
\ea\right.
\eeq{cLaSS}
for all $T>0$. Moreover,
\be
\lim_{t\to\infty}\left(\|\nabla\bfu(t)\|_2+\|\bfu(t)\|_6+|\dot{\bfchi}(t)|+|\bfchi(t)|\right)=0.
\eeq{asympt}
\ET{5.1_01}

Before giving the (lengthy) proof of \theoref{5.1_01}, postponed until Sections \ref{techLem} and \ref{proofstab},
several comments are in order.
\Br
We first point out that the initial datum in \eqref{04_111} is assumed to satisfy $\bfu^0\in W^{1,2}(\Omega)$.
Overall, the statement may appear unsatisfactory since it {\it assumes} that
\be
\lambda_2^{-1}(\lambda)<\infty\qquad\mbox{and}\qquad\lambda<\lambda_2
\eeq{asym}
but, as we now discuss, a stronger result appears {\it in general} out of reach.
It is readily seen that the first assumption of \eqref{asym} is satisfied whenever
\be
\bfu_0(\lambda)\in L^2(\Omega)\, .
\eeq{u0L2}
Actually, by the triangle and H\"older inequalities we get
\begin{eqnarray*}
\frac1{\lambda_2(\lambda)} &=& \sup_{\mbox{\footnotesize $\ \bfu\in\cald^{1,2}(\real^3)$}}\frac{(\bfu\cdot\nabla\bfu,\bfu_0(\lambda))-(\hat{\bfu}\cdot\nabla\bfu,\bfu_0(\lambda))}{\|\nabla\bfu\|_2^2}\\
&\le& \sup_{\mbox{\footnotesize $\ \bfu\in\cald^{1,2}(\real^3)$}}\frac{\|\bfu\|_6\cdot\|\bfu_0(\lambda)\|_3+|\hat{\bfu}|\cdot\|\bfu_0(\lambda)\|_2}{\|\nabla\bfu\|_2}<\infty
\end{eqnarray*}
since $\|\nabla\bfu\|_2$ bounds both $\|\bfu\|_6$ and $|\hat{\bfu}|$ (by \lemmref{1.1}). This proves the first of \eqref{asym} whenever \eqref{u0L2} holds.
However, \theoref{exi} {ensures} {\it in general} only that $\bfu_0(\lambda)\in L^q(\Omega)$ for all $q>2$,
plus some integrability conditions on its derivatives. At the same time, we can readily show that
\be
\exists\,\bfu_0\in L^q(\Omega)\,, \mbox{all $q>2$},\,\bfu_0\not\in L^2(\Omega),\ \exists\,\bfu\in\cald^{1,2}(\real^3)\,\ \mbox{s.t.}\,\
\frac{((\bfu-\hat{\bfu})\cdot\nabla\bfu,\bfu_0)}{\|\nabla\bfu\|_2^2}=\infty\, .
\eeq{u0notL2}
Actually,  take $\bfu\in\cald^{1,2}(\real^3)$ and $\bfu_0$  such that,
as $|x|\to\infty$,
\be	
|\bfu(x)|\asymp\frac{c}{|x|^{1/2}(\ln|x|)^{2/3}}\, ,\quad|\nabla\bfu(x)|\asymp\frac{c}{|x|^{3/2}(\ln|x|)^{2/3}}\, ,\quad
|\bfu_0(x)|\asymp\frac{c}{|x|^{3/2}}\,.
\label{Fi}
\end{equation}
If we split the fraction as
$$
\frac{(\bfu\cdot\nabla\bfu,\bfu_0)-(\hat{\bfu}\cdot\nabla\bfu,\bfu_0)}{\|\nabla\bfu\|_2^2}\, ,
$$
by the H\"older inequality as above, and (\ref{Fi}) we find that
$$
\sup_{\mbox{\footnotesize $\ \bfu\in\cald^{1,2}(\real^3)$}}\frac{(\bfu\cdot\nabla\bfu,\bfu_0)}{\|\nabla\bfu\|_2^2}<\infty\, .
$$
On the other hand,
we also have
$$
\frac{|(\hat{\bfu}\cdot\nabla\bfu,\bfu_0)|}{\|\nabla\bfu\|_2^2}=+\infty\qquad\forall\hat{\bfu}\in\real^3\setminus\{\bf0\}\, .
$$
This proves \eqref{u0notL2}. Incidentally, we notice that the derivatives of the above $\bfu_0$ satisfy
$$
|\nabla\bfu_0(x)|\asymp\frac{c}{|x|^{5/2}}\, ,\quad|D^2\bfu_0(x)|\asymp\frac{c}{|x|^{7/2}}\, ,\quad\mbox{as }|x|\to\infty\, ,
$$
so that \eqref{sfaco} is fulfilled (in fact, we even have larger intervals for $r$ and $s$).
In conclusion, we just saw that \eqref{u0L2} gives a {\it sufficient} condition for the validity of the first of
\eqref{asym}. This condition may not be necessary but the above example suggests that \eqref{asym}$_1$ could fail if \eqref{sfaco2} is evaluated along a {\it generic} solution $\bfu_0$.
\ER{max}
\Br
Once the first condition in \eqref{asym} is satisfied, in order to apply \theoref{5.1_01} one needs to check the second condition.
We already observed that, in general, \eqref{lambda2ge0} holds.
The assumptions of \theoref{5.1_01} require slightly more, namely
$$0<\lambda_2(\lambda)\le \lambda_1(\lambda).$$
This means that if $0<\lambda<\lambda_2(\lambda)$, then the corresponding steady-state solution ${\sf s}_0(\lambda)$
determined in \theoref{exi} is unique {\it and} stable.
\ER{max1}
\subsubsection{Some technical lemmas}\label{techLem}

We prove here some preliminary results. Recall that $\Omega=\real^3\setminus\Omega_0$ and $\Omega_R=\Omega \cap B_R$.
\Bl Let $(\bfu,p)\in [\cald^{1,2}(A)\cap W^{2,2}(D)]\times D^{1,2}(D)$ be such that $\bfu|_{\Omega_0}=\hat{\bfu}$ for some $\hat{\bfu}\in\real^3$ with either $\{A,D\}\equiv\{\real^3,\Omega\}$ or $\{A,D\}\equiv\{B_R,\Omega_R\}$.
Then, there exists a constant $C>0$, depending only on the regularity of $\Omega$, such that
$$
\|D^2\bfu\|_{2,D}+\|\nabla p\|_{2,D}\le C\,(\|\Div \mathbb T(\bfu,p)\|_{2,D}+\|\nabla\bfu\|_{2,D}+|\hat{\bfu}|)\,.
$$
\EL{stoces}
{\it Proof.} See \cite[Lemma 1]{Hey} where the domain is requested to be of class $C^3$. However, $C^2$ suffices.\hfill$\square$\par\smallskip

\Bl Let $\bfu\in \cald^{1,2}(A)\cap W^{2,2}(D)$, with  $\bfu|_{\Omega_0}=\hat{\bfu}$ for some $\hat{\bfu}\in\real^3$,
$A$ and $D$ as in \lemmref{stoces}, and  let $\bfv\in L^2(D)$. Then, for any $\varepsilon>0$ there exists a positive constant $C$ (depending only on $\varepsilon$, $\bfu_0$, and the regularity of $\Omega$) such that

\begin{multline*}
\big|\left(\partial_1\bfu-\bfu_0\cdot\nabla\bfu+(\hat{\bfu}-\bfu)\cdot\nabla(\bfu_0-\bfu),\bfv\right)_D\big|\\
\hspace*{29mm}\le C\,(\|\nabla\bfu\|_{2,D}^2 +\|\nabla\bfu\|_{2,D}^4+\|\nabla\bfu\|_{2,D}^6)+\varepsilon(\|D^2\bfu\|_{2,D}^2+\|\bfv\|_{2,D}^2).
 \end{multline*}

\EL{prep}
{\it Proof.} Let us denote by $I_i$, $i=1,\ldots6$, in the order, the six terms in the scalar product. Taking  into account that, by \theoref{exi}, $\bfu_0\in L^\infty(\Omega)\cap D^{1,q}(\Omega)$, $q=2,3$, and using \eqref{1.9}, \eqref{1.10},  H\"older and Cauchy-Schwarz inequalities we readily get
$$
\sum_{i=1}^5\left|I_i\right|\le C\,( \|\nabla\bfu\|_{2,D}^2 +\|\nabla\bfu\|_{2,D}^4)+\half\varepsilon\|\bfv\|_{2,D}^2\,.
$$
Moreover, again by H\"older inequality, \cite[Lemma 1]{Hey} (see also the ``proof'' of \lemmref{stoces}),
\eqref{1.9} and \remref{2.1},
$$\ba{rl}\medskip
|I_6|\le \|\bfu\|_6\|\nabla\bfu\|_3\|\bfv\|_2&\!\!\!\le C\,\|\nabla\bfu\|_2\left(\|D^2\bfu\|_{2,D}^{\frac12}\|\nabla\bfu\|_{2,D}^{\frac12}+\|\nabla\bfu\|_{2,D}\right)\|\bfv\|_{2,D}\\
&\!\!\!\le C\,(\|\nabla\bfu\|_{2,D}^4+\|\nabla\bfu\|_{2,D}^6)+\varepsilon\,\|D^2\bfu\|_{2,D}^2+\half\varepsilon\,\|\bfv\|_{2,D}^2\,.
\ea
$$
The lemma is proved.\hfill$\square$\par\smallskip

We will also need the following technical result.

\Bl
Let $y:[0,\infty)\mapsto[0,\infty)$, be absolutely continuous and satisfying
\be
\dot{y}(t)\le a(t)+b(t)[y(t)+y^\alpha(t)]\,,\ \ \alpha\ge 1,\ \ \mbox{for a.e. }t>0\,,
\eeq{a}
where $a,b\in L^\infty(0,\infty)$.
Assume, also, $y\in L^1(0,\infty)$, and set
$$
{\sf a}:=\essup{t\in (0,\infty)}\,|a(t)|\,,\quad {\sf b}:=\essup{t\in (0,\infty)}\,|b(t)|\,.
$$
Then, there exists $\delta>0$, such that if
\be
y(0)\le \delta\,,\ \ \int_0^\infty y(s)\,ds\le\delta^2
\eeq{b}
it follows that:
$$\mbox{$y(t)< M\,\delta$ for all $t\in (0,\infty)\,,$ $M:=3\max\{1,2{\sf a},2{\sf b}\}$}\,.
$$
\EL{gronwa}
{\it Proof.} Let
$$
Y:=y^2\,,\ \beta:=(1+\alpha)/2\,.
$$
Multiplying both sides of \eqref{a} by $y$  we get
\be
\dot{Y}\le 2{\sf a}y+2{\sf b}[Y+Y^\beta].
\eeq{a1}
Contradicting the statement means that there exists $t_0>0$ such that \be y(0)\le\delta\,,\,\ y(t)<M\,\delta\,, \ \ \mbox{for all $t\in (0,t_0)$\,,\, and  \,$y(t_0)=M\,\delta.$}\eeq{c} Then, integrating both sides of \eqref{a1} from 0 to $t_0$ and using the latter and \eqref{b}$_2$, we deduce, in particular
$$\ba{rl}\medskip
Y(t_0)&\!\!\!\le Y(0) +2{\sf a}\Int0{\infty}y(s)\,{\rm d}s+2{\sf b}\Int0{t_0} Y(s)\,{\rm d}s +2{\sf b}\Int0{t_0}Y^\beta(s)\,{\rm d}s\,\\
&\!\!\!\le \delta^2+2{\sf a}\delta^2+2{\sf b}M\delta^3
+2{\sf b}\delta^2(M\delta)^\alpha\le\frac{M}3 \delta^2 \,(2+M\delta+(M\delta)^\alpha)\,.
\ea
$$
Therefore, choosing $\delta>0$ in such a way that
$$
2+M\delta+(M\delta)^\alpha<3M
$$
we deduce $y(t_0)<M\delta$, which contradicts \eqref{c}$_3$. \hfill$\square$\par\smallskip

\subsubsection{Proof of {\theoref{5.1_01}}}\label{proofstab}
\

{\bf Part 1: existence.}
To prove the existence of a solution to \eqref{04_111},
we follow the arguments introduced and developed in \cite{GaSi1,GaSi2}. Let $\{\Omega_R\}$, $R\in\nat$, be an increasing sequence such that $\Omega=\cup_{R\in\nat}\Omega_R$ and, for each fixed $R$, consider the  problem
\be\ba{cc}\medskip\left.\ba{lr}\medskip
& \partial_t\bfu_R-\Div\mathbb T(\bfu_R,{p}_R)\\
&=\lambda\Big[\partial_1\bfu_R-\bfu_0\cdot\nabla\bfu_R+({\bfsigma}_R-\bfu_R)\cdot\nabla(\bfu_0-\bfu_R)\Big]\\
& \Div\bfu_R=0\ea\right\}\, \ \mbox{in $\Omega_R\times(0,\infty)$}\,,\\ \medskip
\bfu_R|_{\partial\Omega}={{\bfsigma}_R}(t)\,, \ \ \bfu_R|_{\partial B_R}=\0\,, \ \ \mbox{in $(0,\infty)$}\\ \medskip
\dot{\bfsigma}_R+\omega_{\sf n}^2\bfchi_R+\varpi\Int{\partial\Omega}{} \mathbb T(\bfu_R,{p}_R)\cdot\bfn=\0\,, \ \ \dot{\bfchi}_R=\bfsigma_R\ \ \mbox{in $(0,\infty)$}\\
\bfu_R(x,0)=\bfu^0\,,\ x\in\Omega_R\,,\ \ \bfchi_R(0)=\bfchi^{0}\,,\ \ {\bfsigma}_R(0)=\bfchi^1\,.
\ea
\eeq{04_11R}

Our approach to existence develops in two steps. In the first step, by the classical Galerkin method we show that \eqref{04_11R} has a solution in the class \eqref{cLaSS}. This is accomplished with the help of a suitable base, constituted by eigenvectors of a modified Stokes problem. This procedure also leads to the proof of estimates for $(\bfu_R,p_R,\bfchi_R)$ with bounds that are independent of $R$. In this way, in the second step, we will pass to the limit $R\to\infty$ and show that the limit functions $(\bfu,p,\bfchi)$ solve the original problem \eqref{04_111} along with the asymptotic property \eqref{asympt}.

We start putting \eqref{04_11R} in a ``weak'' form. If we multiply \eqref{04_11R}$_1$ by $\bfpsi\in\cald^{1,2}(B_R)$, integrate by parts and use \eqref{04_11R}$_{2,3,4}$, we deduce
\be\ba{ll}\medskip
(\partial_t\bfu,\bfpsi)+2(\mathbb D(\bfu),\mathbb D ({\bfpsi}))+\varpi^{-1}(\dot{\bfsigma}+\omega^2_{\sf n}\bfchi)\cdot{\hat{\bfpsi}}\\
=\lambda\,\left[\partial_1\bfu-\bfu_0\cdot\nabla\bfu+({\bfsigma}-\bfu)\cdot\nabla(\bfu_0 -\bfu),\bfpsi\right]\ \ \mbox{for all $\bfpsi\in\cald^{1,2}(B_R)$}\,,
\ea
\eeq{wefo}
where $\bfsigma=\dot{\bfchi}$ and, for simplicity, the subscript ``$R$'' has been omitted and $(\cdot,\cdot)\equiv(\cdot,\cdot)_{\Omega_R}$. Using  standard procedures \cite{Gah}, one finds that if $(\bfu,p,\bfchi,\bfsigma)$ is  a smooth solution to \eqref{wefo}, then  it also satisfies \eqref{04_11R}$_{1-6}$.
In \cite{GaSi1} it is shown that the problem
\be
\begin{array}{l}
\left.
\begin{array}{l}
\hspace{-0.2cm} \displaystyle  - \nabla \cdot \mathbb T(\bfpsi,\phi)=\mu  \,\bfpsi \medskip \\
\hspace{-0.2cm} \displaystyle \Div   \bfpsi = 0 \medskip
\end{array}
\right\} \mbox{ in } \Omega_R \,,\\  \medskip
\displaystyle \bfpsi  = \hat{\bfpsi} \ \   \mbox{ on } \partial\Omega\,, \ \
\displaystyle \bfpsi= 0\ \   \mbox{ on } \partial B_R \,, \\
\displaystyle \mu\,  \hat{\bfpsi}=\varpi\int_{\partial\Omega} \mathbb T(\bfpsi,\phi)\cdot \bfn \,,
\end{array}
\eeq{eipr}
with the natural extension $\bfpsi(x)=\hat{\bfpsi}$ in $\Omega_0$,
admits a denumerable number of positive eigenvalues $\{\mu_{Ri}\}$ clustering at infinity,
and  corresponding eigenfunctions $\{\bfpsi_{Ri}\}_i \subset {\mathcal D}^{1,2}(B_R)\cap W^{2,2}(\Omega_R)$    forming an orthonormal basis of ${\mathcal{H}}(B_R)$ that is also orthogonal in $\cald^{1,2}(\Omega_R)$. Also, the correspondent ``pressures'' satisfy $\phi_{Ri}\in W^{1,2}(\Omega_R)$, $i\in\nat$. Thus, for each fixed $R\in\nat$, we look for ``approximated'' solutions to \eqref{wefo} of the form
\be
\bfu_N(x,t)=\sum_{k=1}^Nc_{kN}(t)\bfpsi_{Rk}(x)\,,\ \ \bfsigma_N(t)=\sum_{k=1}^Nc_{kN}(t)\hat{\bfpsi}_{Rk}\,,\ \ \bfchi_N(t)\,,
\eeq{gal0}
where the vector functions $\bfc_N(t):=\{c_{1N}(t),\ldots c_{NN}(t)\}$ and $\bfchi_N(t)$ satisfy the following system of  equations $(i=1,\ldots,N)$
\be
\begin{split}
(\partial_t\bfu_N,\bfpsi_{Ri})+2(\mathbb D(\bfu_N),\mathbb D ({\bfpsi}_{Ri}))+
\varpi^{-1}(\dot{\bfsigma}_N+\omega^2_{\sf n}\bfchi_N)\Cdot{\hat{\bfpsi}}_{Ri}\\
=\lambda\,\left[\partial_1\bfu_N-\bfu_0\Cdot\nabla\bfu_N+({\bfsigma}_N-\bfu_N)\Cdot\nabla(\bfu_0-\bfu_N),\bfpsi_{Ri}\right]\, ,
\end{split}
\eeq{gal}
with $\bfsigma_N=\dot{\bfchi}_N$.
Indeed, \eqref{gal} is a system of first order differential equations in normal form in the unknowns $\bfc_N,\bfchi_N$. To this end,
it suffices to observe that
\be
\langle\bfpsi_{Ri},\bfpsi_{Rj}\rangle:=(\bfpsi_{Ri},\bfpsi_{Rj})+\varpi^{-1}\hat{\bfpsi}_{Ri}\cdot\hat{\bfpsi}_{Rj}=\delta_{ij}\,,
\eeq{ortg}
so that the derivatives with respect to time can be grouped and \eqref{gal} is equivalent to the system
\be\ba{l}
\dot{c}_{iN}= F_i(\bfc_N,\bfchi_N)\,,\ \ i=1,\ldots,N\,,\\ \medskip
F_i:=\displaystyle\sum_{k=1}^Nc_{kN}\Big[\lambda\big(\partial_1\bfpsi_{Rk}-\bfu_0\cdot\nabla\bfpsi_{Rk}+(\hat{\bfpsi}_{Rk}-\bfpsi_{Rk})\cdot\nabla\bfu_0,\bfpsi_{Ri}\big) \\ \medskip
-2\big(\mathbb D(\bfpsi_{Rk}),\mathbb D(\bfpsi_{Ri})\big)\Big]
-\lambda\displaystyle\sum_{k,m=1}^Nc_{kN}c_{mN}\left((\hat{\bfpsi}_{Rk}-{\bfpsi}_{Rk})\cdot\nabla{\bfpsi}_{Rm},{\bfpsi}_{Ri}\right)\\ -\frac{\omega_{\sf n}^2}{\varpi}\bfchi_N\cdot\hat{\bfpsi}_{Ri}\,,
\ea
\eeq{ODE}
which we equip with the following initial conditions:
\be
c_{iN}(0)=(\bfu^0,\bfpsi_{Ri})+\varpi^{-1}\bfchi^1\cdot\hat{\bfpsi}_{Ri}\,,\ \ \bfchi_N(0)=\bfchi^0\,.
\eeq{inco}
From \eqref{gal0} and \eqref{ortg}, it follows that
\be
\|\bfu_N(0)\|_{2,\Omega_R}^2+\varpi^{-1}|\bfsigma_N(0)|^2\le\|\bfu^0\|_{2,\Omega}^2+\varpi^{-1}|\bfchi^1|^2\,.
\eeq{incou}
Likewise, since
$$
2 (\mathbb D(\bfpsi_{Ri}),\mathbb D(\bfpsi_{Rj}))= \mu_{Ri} \Big[ \varpi^{-1} \hat{\bfpsi}_{Ri} \cdot \hat{\bfpsi}_{Rj} +(\bfpsi_{Ri},\bfpsi_{Rj})\Big] = \mu_{Ri} \delta_{ij}
$$
we have
\begin{eqnarray*}
\mathbb D(\bfu_N(0))&= &\sum_{j=1}^{N}c_{jN}(0)\mathbb D(\bfpsi_{Rj}) = 2 \sum_{j=1}^{k}  \frac{1}{\mu_{Rj}} (\mathbb D(\bfu^0), \mathbb D(\bfpsi_{Rj}))_{\Omega_R} \mathbb  D(\bfpsi_{Rj})\\ &=& \sum_{j=1}^{N}   \frac{(\mathbb D(\bfu^0), \mathbb D(\bfpsi_{Rj}))_{\Omega_R}}{\| \mathbb D(\bfpsi_{Rj}) \|_{2,\Omega_R}^{2}} \,   \mathbb D(\bfpsi_{Rj})
\end{eqnarray*}
and, therefore,
\be
\|\mathbb D( \bfu_N(0))\|_{2,\Omega_R}\leq \|\mathbb D(\bfu^{0})\|_{2,\Omega}.
\eeq{incod}

We shall now derive three basic ``energy estimates''.
Multiplying both sides of \eqref{gal}$_1$ by $c_{iN}$, summing over $i$,  integrating by parts over $\Omega_R$ and using \eqref{gal0} and \eqref{1.8}, we get
$$
\half\ode{}t\left[\|\bfu_N\|_2^2+\varpi^{-1}(|\bfsigma_N|^2+\omega_{\sf n}^2|\bfchi_N|^2)\right]
+\|\nabla(\bfu_N)\|_2^2-\lambda\,((\bfu_N-\bfsigma_N)\cdot\nabla\bfu_N,\bfu_0) = 0\,,
$$
which, by \eqref{04_111}, \remref{2.1} and \eqref{sfaco2}, in turn furnishes
\be
\half\ode{}t\left[\|\bfu_N\|_2^2+\varpi^{-1}(|\bfsigma_N|^2+\omega_{\sf n}^2|\bfchi_N|^2)\right]+ \gamma \,\|\nabla\bfu_N\|_2^2\le 0,
\eeq{gal1}
where we have set
\be
\gamma :=1-\frac{\lambda}{\lambda_2}>0.
\eeq{gamma}

We next multiply both sides of \eqref{gal}$_1$ by $\dot{c}_{iN}$, sum over $i$,  and integrate by parts over $\Omega_R$ as necessary.
Taking again into account \eqref{gal0}, \eqref{1.8} and \remref{2.1}, we obtain
\be
\begin{split}
\half\ode{}t\|\nabla\bfu_N\|_2^2+\|\partial_t\bfu_N\|_2^2+\varpi^{-1}|\dot{\bfsigma}_N|^2
=\lambda\,\Big(\partial_1\bfu_N-\bfu_0\cdot\nabla\bfu_N \\ +({\bfsigma}_N-\bfu_N)\cdot\nabla\bfu_0   -({\bfsigma}_N-\bfu_N)\cdot\nabla\bfu_N,\partial_t\bfu_N\Big)-\frac{\omega^2_{\sf n}}{\varpi}\bfchi_N\cdot\dot{\bfsigma}_N\,. &
\end{split}
\eeq{gal2}

Finally, we multiply both sides of \eqref{gal}$_1$ by $\lambda_{Ri}c_{iN}$ and sum over $i$. Integrating by parts over $\Omega_R$ and employing \eqref{eipr} and, one more time, \eqref{04_111}, \eqref{1.8}, and  \remref{2.1} we show
\be
\begin{split}
\half\ode{}t\|\nabla \bfu_N\|_2^2+\|\Div\mathbb T(\bfu_N,p_N)\|_2^2+\varpi\,|\bfS_N|^2
=\lambda\,\Big(\partial_1\bfu_N-\bfu_0\cdot\nabla\bfu_N\\ +({\bfsigma}_N-\bfu_N)\cdot\nabla\bfu_0 -({\bfsigma}_N-\bfu_N)\cdot\nabla\bfu_N,\Div\mathbb T(\bfu_N,p_N)\Big)
+{\omega^2_{\sf n}}\bfchi_N\cdot\bfS_N
\end{split}
\eeq{gal3}
where
$$
\bfS_N:=\int_{\partial\Omega}\mathbb T(\bfu_N,p_N)\cdot\bfn\,,\ \ p_N:=\sum_{k=1}^Nc_{kN}\phi_{Rk}\,.
$$

We shall now derive a number of estimates for the approximated solutions,
paying attention that the constants involved are independent of $N$ and $R$. Such generic constants will be denoted by $C$, which can thus depend, at most, on $\Omega$, $\bfu_0$ and the physical constants involved in the problem. Moreover, without specification, its value may change from a line to the next one (e.g.\ $2C\le C$). From \eqref{gal1}, \eqref{incou} and \eqref{1.10} we get
\be\ba{ll}
\Sup{t\in(0,\infty)}\left[\|\bfu_N(t)\|_2^2+\varpi^{-1}(|\bfsigma_N(t)|^2+\omega_{\sf n}^2|\bfchi_N(t)|^2)\right]\\ +\gamma\Int0\infty(2\kappa_1^{-2}|\bfsigma_N(s)|^2+\|\nabla\bfu_N(s)\|_2^2)\,{\rm d}s \\ \le\|\bfu^0\|_2^2+\varpi^{-1}(|\bfchi^1|^2+\omega_{\sf n}^2|\bfchi^0|^2)\,.
\ea
\eeq{est1}
Such an estimate implies, in particular, that the initial-value problem \eqref{ODE}--\eqref{inco} has a (unique) solution in the whole interval $(0,\infty)$.
Moreover, from \eqref{gal2}, Cauchy-Schwarz inequality and \lemmref{prep} with $\bfv\equiv\partial_t\bfu_N$ and $\varepsilon\equiv\varepsilon_1<\half$, we infer
\arraycolsep=1pt
\be\ba{rl}
\ode{}t\|\nabla\bfu_N\|_2^2+\half\|\partial_t\bfu_N\|_2^2+|\dot{\bfsigma}_N|^2 \le & C\left(\|\nabla\bfu_N\|_2^2+\|\nabla\bfu_N\|_2^4+\|\nabla\bfu_N\|_2^6+|\bfchi_N|^2\right)\\ \medskip & +\varepsilon_1\,\|D^2\bfu_N\|_2^2\\ \medskip
\le & C\left(\|\nabla\bfu_N\|_2^2+\|\nabla\bfu_N\|_2^6+|\bfchi_N|^2\right)+\varepsilon_1\,\|D^2\bfu_N\|_2^2.
\ea
\eeq{est2}
Likewise, employing \lemmref{prep}, this time with $\bfv\equiv\Div \mathbb T(\bfu_N,p_N)$, from \eqref{gal3} we obtain
$$\ba{rl}\medskip
\half \ode{}t\|\nabla\bfu_N\|_2^2+\|\Div\mathbb T(\bfu_N,p_N)\|_2^2+\varpi\,|\bfS_N|^2
\le & C\left(\|\nabla\bfu_N\|_2^2+\|\nabla\bfu_N\|_2^6+|\bfchi_N|^2\right)\\ & +\varepsilon_2\,\|D^2\bfu_N\|_2^2
\ea
$$
\arraycolsep=5pt
which, in turn, combined with \lemmref{stoces} implies, by taking $\varepsilon_2$ small enough
\be
\half\ode{}t\|\nabla\bfu_N\|_2^2+C\,\|D^2\bfu_N\|_2^2+\varpi\,|\bfS_N|^2
\le C\left(\|\nabla\bfu_N\|_2^2+\|\nabla\bfu_N\|_2^6+|\bfchi_N|^2\right)\,.
\eeq{est3}
Summing side-by-side \eqref{est2} and \eqref{est3} by choosing $\varepsilon_1$ sufficiently small we deduce
\be
\begin{split}
\ode{}t\|\nabla\bfu_N\|_2^2+C\,(\|D^2\bfu_N\|_2^2+\|\partial_t\bfu_N\|_2^2+|\dot{\bfsigma}_N|^2+\varpi\,|\bfS_N|^2)\\
\le C\left(\|\nabla\bfu_N\|_2^2+\|\nabla\bfu_N\|_2^6+|\bfchi_N|^2\right),
\end{split}
\eeq{est4}
which furnishes,
in particular,
\be
\ode{}t\|\nabla\bfu_N\|_2^2
\le C\left(\|\nabla\bfu_N\|_2^2+\|\nabla\bfu_N\|_2^6+|\bfchi_N|^2\right).
\eeq{est5}

In \eqref{est5} we put $y(t)=\|\nabla\bfu_N(t)\|_2^2$, $a(t)=C\,|\bfchi_N(t)|^2$ and $b=C$. By \eqref{est1}, it follows that both $a$ and $b$ satisfy the assumptions of \lemmref{gronwa} with $\alpha=3$. Hence, there exists $\delta>0$ such that if \eqref{b} holds, namely
\be
\|\nabla\bfu_N(0)\|_{2,\Omega_R}^2\le \delta\,,\ \ \int_0^\infty\|\nabla\bfu_N(s)\|_{2,\Omega_R}^2\,ds\le\delta^2\, ,
\eeq{b11}
then
\be
\sup_{t\in(0,\infty)}\|\nabla\bfu_N(t)\|_{2,\Omega_R}\le M\delta\,.
\eeq{fiestN}
We take $\varepsilon>0$ in \eqref{smalleps} such that
$$
\varepsilon^2\le\min\Big\{\delta,\gamma\delta^2\min\{1,\varpi(1+\omega_{\sf n}^2)^{-1}\}\Big\}\, .
$$
Then, from \eqref{incod} and \eqref{smalleps} we know that
$$
\|\nabla\bfu_N(0)\|_{2,\Omega_R}^2\le\|\nabla\bfu^0\|_{2,\Omega}^2\le\big(\|\bfu^0\|_{1,2}+|\bfchi^0|+|\bfchi^1|\big)^2\le\varepsilon^2\le\delta\, ,
$$
while from \eqref{est1} we infer that
\begin{eqnarray*}
\Int0\infty \|\nabla\bfu_N(s)\|_{2,\Omega_R}^2\,{\rm d}s &\le& \gamma^{-1}\Big[\|\bfu^0\|_2^2+\varpi^{-1}(|\bfchi^1|^2+\omega_{\sf n}^2|\bfchi^0|^2)\Big]\\
&\le& \gamma^{-1}\max\left\{1,\frac{1+\omega_{\sf n}^2}{\varpi}\right\}\Big[\|\bfu^0\|_2^2+|\bfchi^1|^2+|\bfchi^0|^2\Big]\\
&\le& \gamma^{-1}\max\left\{1,\frac{1+\omega_{\sf n}^2}{\varpi}\right\}\, \varepsilon^2\le\delta^2\, .
\end{eqnarray*}

Therefore, with the above choice of $\varepsilon>0$, both conditions in \eqref{b11} are satisfied and \eqref{fiestN} holds with a constant
$C_0=M\delta>0$, independent of $N$ and $R$, namely
\be
\sup_{t\in(0,\infty)}\|\nabla\bfu_N(t)\|_2\le C_0\,.
\eeq{fiest}

Employing \eqref{fiest} in \eqref{est4} and keeping in mind \eqref{est1} we conclude
\be
\int_0^T\left(\|D^2\bfu_N(s)\|_2^2+\|\partial_t\bfu_N(s)\|_2^2\right)\,{\rm d}s\le C_1\,T\,,\ \ \mbox{for all $T>0$}\,,
\eeq{scest}
with $C_1$ another positive constant independent of $N$ and $R$. Thanks to \eqref{est1} and \eqref{scest}, we can now use a standard argument (see e.g. \cite{GaSi1}) to prove the existence of a subsequence $\{(\bfu_{N_k},\bfchi_{N_k},\bfsigma_{N_k})\}$ converging in suitable topology to some $(\bfu_{R},\bfchi_{R},\bfsigma_{R})$ in the class \eqref{cLaSS} (with $\Omega$ replaced by $\Omega_R$ and $\real^3$ replaced by $B_R$)  and satisfying \eqref{wefo}. Since, clearly, $(\bfu_{R},\bfchi_{R},\bfsigma_{R})$ continues to obey the bounds \eqref{est1} and \eqref{scest}, we can similarly select a subsequence $(\bfu_{R_m},\bfchi_{R_m},\bfsigma_{R_m})$ converging (again, in suitable topology) to a certain $(\bfu,\bfchi,\bfsigma)$ that is in the class \eqref{cLaSS} and obeys \eqref{est1}, \eqref{scest}, and \eqref{04_111} for a.e. $x\in \Omega$ and $t\in(0,\infty)$. The demonstration of this convergence is rather typical and we omit it, referring to \cite[Step 3 at p. 141]{GaSi1} for details. Thus, the proof of existence is completed.

{\bf Part 2: uniqueness.}  This part of the proof is quite standard and we only sketch it here. Let $(\bfu_i,p_i,\bfchi_i)$, $i=1,2$, be two solutions to \eqref{04_111} corresponding to the same initial data, and set $\bfu:=\bfu_1-\bfu_2$, $p=p_1-p_2$, $\bfchi=\bfchi_1-\bfchi_2$. We thus have
\be\ba{cc}\left.
\ba{rr}
& \partial_t\bfu-\Delta\bfu+\nabla {p}=\lambda\,\Big[\partial_1\bfu-\bfu_0\cdot\nabla\bfu\\
& +(\dot{\bfchi}-\bfu)\cdot\nabla(\bfu_0-\bfu_2) -(\dot{\bfchi}_1  -\bfu_1)\cdot\nabla\bfu\Big]\\ \medskip
& \Div\bfu=0\ea\right\}\ \ \mbox{in $\Omega\times(0,\infty)$}\,,\\ \medskip
\bfu(x,t)={\dot{\bfchi}}(t)\,, \ \mbox{ $(x,t)\in\partial\Omega\times(0,\infty)$}\,,\\ \medskip
\ddot{\bfchi}+\omega_{\sf n}^2\bfchi+\varpi\Int{\partial\Omega}{} \mathbb T(\bfu,{p})\cdot\bfn=\0\,, \ \ \mbox{in $(0,\infty)$}\\
\bfu(x,0)=\0\,,\ x\in\Omega\,,\ \ \bfchi(0)=\0\,,\ \ \dot{\bfchi}(0)=\0
\,.
\ea
\eeq{04_11U}
We dot-multiply \eqref{04_11U}$_1$ by $\bfu$, integrate by parts over $\Omega$ and use \eqref{04_11U}$_{2-4}$ to obtain
\be
\half\ode{}t\left[\|\bfu\|_2^2+\varpi^{-1}(|\dot{\bfchi}|^2+\omega_{\sf n}^2|\bfchi|^2)\right]+\gamma \|\nabla\bfu\|_2^2\le \lambda\,((\dot{\bfchi}-\bfu)\cdot\nabla\bfu_2,\bfu)\,,
\eeq{un1}
where we recall that $\gamma$ is defined in \eqref{gamma}. From \eqref{1.10} and Cauchy--Schwarz inequality we get
$$
|(\bfchi\cdot\nabla\bfu_2,\bfu)|\le \half\gamma\|\nabla\bfu\|_2^2+c\,\|\nabla\bfu_2\|_2^2\|\bfu\|_2^2\,,
$$
whereas from \eqref{1.9}, H\"older,  Sobolev and Cauchy--Schwarz inequalities,
$$
|(\bfu\cdot\nabla\bfu_2,\bfu)|\le \|\bfu\|_6\|\nabla\bfu_2\|_3\|\bfu\|_2\le \half\gamma \|\nabla\bfu\|_2^2+c\,\|\bfu_2\|^2_{2,2}\|\bfu\|_2^2\,.
$$
Replacing the last two displayed relations in \eqref{un1}, we thus conclude
\be
\ode Et\le c\, g(t)\, E(t)
\eeq{EN}
where $g:=\|\bfu_2\|_{2,2}^2$, $E:=
\|\bfu\|_2^2+\varpi^{-1}(|\dot{\bfchi}|^2+\omega_{\sf n}^2|\bfchi|^2)$. Since  $\bfu_2$ is in the class \eqref{cLaSS}, we have $g\in L^1(0,T)$, for all $T>0$ and also, by assumption, $E(0)=0$. Uniqueness then follows by using Gronwall's lemma in \eqref{EN}.

{\bf Part 3: stability.} We finally prove the validity of \eqref{asympt}. In this regard, we begin to observe that the solution just constructed satisfies, in particular,
\be
\Sup{t\in(0,\infty)}\left(\|\bfu(t)\|_2+\|\nabla\bfu(t)\|_2+|\bfchi(t)|+|\dot{\bfchi}(t)|\right)+\Int0\infty(|\dot{\bfchi}(s)|^2+\|\nabla\bfu(s)\|_2^2)\,{\rm d}s\le K\,,
\eeq{asy1}
where $K>0$ is a constant depending only on the data. By dot-multiplying both sides of \eqref{04_111}$_1$ by $\partial_t\bfu$ and proceeding as in the proof of \eqref{est2} we obtain
$$\ba{ll}\medskip
\ode{}t(\|\nabla\bfu\|_2^2+\frac{\omega^2_{\sf n}}{\varpi}\bfchi\Cdot\dot{\bfchi})+\|\partial_t\bfu\|_2^2+\varpi^{-1}|\ddot{\bfchi}|^2\\\hspace*{2cm}=\lambda\,\left(\partial_1\bfu-\bfu_0\Cdot\nabla\bfu+(\dot{\bfchi}-\bfu)\Cdot\nabla\bfu_0 -({\bfsigma}-\bfu)\Cdot\nabla\bfu,\partial_t\bfu\right)+\mbox{$\frac{\omega^2_{\sf n}}{\varpi}$}|\dot{\bfchi}|^2\,.
\ea
$$
We now use, on the right-hand side of this relation, \lemmref{prep} with $\bfv\equiv\partial_t\bfu$, $\varepsilon\equiv\varepsilon_1<\half$,  along with the uniform bound on $\|\nabla\bfu\|_2$ in \eqref{asy1} to get
\be
\ode{}t(\|\nabla\bfu\|_2^2+\mbox{$\frac{\omega^2_{\sf n}}{\varpi}$}\bfchi\cdot\dot{\bfchi})+\half\|\partial_t\bfu\|_2^2+\varpi^{-1}|\ddot{\bfchi}|^2\le C\,\|\nabla\bfu\|_2^2+\varepsilon_1\,\|D^2\bfu\|_2^2+\mbox{$\frac{\omega^2_{\sf n}}{\varpi}$}|\dot{\bfchi}|^2
\eeq{asy2}
Finally, we test both sides of \eqref{04_111}$_1$ by $-\Div\mathbb T(\bfu,p)$  and apply Cauchy-Schwarz inequality to deduce
\[
\|\Div\mathbb T(\bfu,p)\|_2^2\le
- 2\lambda\,\left(\partial_1\bfu-\bfu_0\cdot\nabla\bfu+(\dot{\bfchi}-\bfu)\cdot\nabla(\bfu_0-\bfu),\Div\mathbb T(\bfu,p)\right) +\|\partial_t\bfu\|_2^2\,.
\]
Employing in this inequality \lemmref{prep} with $\bfv\equiv\Div\mathbb T(\bfu,p)$ along with the bound \eqref{asy1} on $\|\nabla\bfu\|_2$,  we infer
$$
\|\Div\mathbb T(\bfu,p)\|_2^2\le C\,\|\nabla\bfu\|_2^2 +\varepsilon_2\,\|D^2\bfu\|_2^2+\|\partial_t\bfu\|_2^2\,,
$$
which, in turn, with the help of \lemmref{stoces} and by selecting $\varepsilon_2$ small enough, entails
\be
\|D^2\bfu\|_2^2+\|\nabla p\|_2^2\le C\,(\|\nabla\bfu\|_2^2 +\|\partial_t\bfu\|_2^2++|\dot{\bfchi}|^2)\,.
\eeq{asy3}
Next, we utilize \eqref{asy3} on the right-hand side of \eqref{asy2} and pick $\varepsilon_2$ suitably, which enables us to find
\be
\ode{}t(\|\nabla\bfu\|_2^2+\mbox{$\frac{\omega^2_{\sf n}}{\varpi}$}\bfchi\cdot\dot{\bfchi})+\mbox{$\frac14$}\|\partial_t\bfu\|_2^2+\varpi^{-1}|\ddot{\bfchi}|^2\le C\,(\|\nabla\bfu\|_2^2+|\dot{\bfchi}|^2)\,.
\eeq{asy4}
Integrating over time both sides of \eqref{asy4}, and taking into account \eqref{asy1}, it follows  that
\be
\partial_t\bfu\in L^2(0,\infty; L^2(\Omega))\,,\ \ \ddot{\bfchi}\in L^2(0,\infty)\,,
\eeq{asy5}
which once replaced in \eqref{asy4}, again with the help of \eqref{asy1},  furnishes
\be
D^2\bfu,\nabla p\in L^2(0,\infty; L^2(\Omega))\,.
\eeq{asy6}
By possibly adding a suitable function of time to $p$, we may get \cite[Theorem II.6.1]{Gab}
\be
p\in L^6(\Omega)\ \ \mbox{and}\,\ \|p(t)\|_6\le C\,\|\nabla p(t)\|_2\,,\ \ \mbox{a.a. $t>0$\,.}
\eeq{p_p}
On the other hand, from \eqref{04_111}$_4$ and standard trace theorems, we have
\[
\begin{split}
|\bfchi(t)|^2 & \le \half\,\left(|\ddot{\bfchi}(t)|^2+\left|\int_{\partial\Omega}\mathbb T(\bfu,p)\cdot\bfn\right|^2\right)\\ & \le C\left(|\ddot{\bfchi}(t)|^2+\|\nabla\bfu\|_{2,2,\Omega_\rho}^2+\|p\|^2_{1,2,\Omega_\rho}\right)\,,
\end{split}
\]
for some fixed $\rho$ which, in vew of \eqref{asy5}--\eqref{p_p}, allows us to conclude that
\be
\bfchi\in L^2(0,\infty)\,.
\eeq{asy7}
Combining \eqref{asy1}, \eqref{asy5} and \eqref{asy7} we get at once
\be
\lim_{t\to\infty}\left(|\bfchi(t)|+|\dot{\bfchi}(t)|\right)=0\,.
\eeq{lim1}
From \eqref{asy1} it follows that there exists at least one unbounded sequence $\{t_n\}\in (0,\infty)$ such that
\be
\lim_{n\to\infty}\|\nabla\bfu(t_n)\|_2=0\,.
\eeq{lim2}
Thus, integrating both sides of \eqref{asy4} between $t_n$ and $t>t_n$ we infer, in particular
\[
\|\nabla\bfu(t)\|_2^2\le  C\left(|\bfchi(t)|\,|\dot{\bfchi}(t)|+|\bfchi(t_n)|\,|\dot{\bfchi}(t_n)|+\int_{t_n}^\infty(\|\nabla\bfu(s)\|_2^2+|\dot{\bfchi}(s)|^2)\,{\rm d}s\right)+ \|\nabla\bfu(t_n)\|_2^2
\]
which, by \eqref{asy1}, \eqref{lim1} and \eqref{lim2} entails
$$
\lim_{t\to\infty}\|\nabla\bfu(t)\|_2=0\,.
$$
The latter and \eqref{1.9} complete the proof of \eqref{asympt}. The proof of \theoref{5.1_01} is completed.

\setcounter{equation}{0}
\subsection{Absence of Oscillatory Solutions below the Stability Threshold}\label{sec:nobif}
As remarked in the previous subsection, if $\lambda_2>0$ and $\lambda<\lambda_2$, the steady solution of \theoref{exi} is unique {\it and} stable.
The objective of this subsection is to show, in addition, that, if $\lambda_2>0$, no oscillatory motion can stem out of the steady-state branch in a suitable function class of solutions $\Gamma$ as long as $\lambda<\lambda_2$. As a direct consequence, a time-periodic bifurcation may occur only at some $\lambda_o>\lambda_2$. More precisely, let ${\sf s}_0(\lambda)$ be the steady-state solution given in \eqref{sfaco}. A generic $T$-periodic solution to \eqref{02} can then always be written as
$$
\bfu(x,t)+\bfu_0(x)\,,\ \ p(x,t)+p_0(x)\,,\ \ \bfchi(t)+\bfchi_0\,,
$$
where $(\bfu,p,\bfchi)$, after the scaling  $\tau=\frac{2\pi}{T}.t$, is a $2\pi$-periodic solution to the following equations
\be\ba{cc}\medskip\left.
\ba{lr}
& \zeta\partial_\tau\bfu-\Delta\bfu+\nabla {p} \\
&= \lambda\,\Big[\partial_1\bfu-\bfu_0\cdot\nabla\bfu+(\zeta\dot{\bfchi}-\bfu)\cdot(\nabla\bfu+\nabla\bfu_0) \Big]\\
& \Div\bfu=0\ea\right\}\ \ \mbox{in $\Omega\times(0,2\pi)$}\,,\\ \medskip
\bfu(x,t)={\zeta\dot{\bfchi}}(t)\,, \ \mbox{ $(x,t)\in\partial\Omega\times(0,2\pi)$}\,,\\ \medskip
\zeta^2\ddot{\bfchi}+\omega_{\sf n}^2\bfchi+\varpi\Int{\partial\Omega}{} \mathbb T(\bfu,{p})\cdot\bfn=\0\,, \ \ \mbox{in $(0,2\pi)$}
\,,
\ea
\eeq{per0}
where $\zeta:=2\pi/T$.
We now introduce the class
\bestar
\begin{split}
\Gamma:= \{(\bfu= &\bar{\bfu}+\bfw,p = \bar{p}+{\sfp}, \bfchi=\bar{\bfchi}+\bfxi): \\
& \bar{\bfu}\in X(\Omega)\,,  \bfw\in\W;\,\ \bar{p}\in W^{1,2}(\Omega)\,, \ {\sfp}\in {\sf P}_\sharp^{1,2} ;\,\ \bfxi\in W^2_\sharp\}\,,
\end{split}
\eestar
which constitutes the functional framework  where, later on, we shall prove the occurrence of a time-periodic bifurcation.
We recall that $\overline{\ \cdot\ }$ denotes the mean value, as defined in \eqref{media}.
In particular, if $\bfu= \bar{\bfu}+\bfw\in \Gamma$, then
$$\overline{\bfw}=\overline{\partial_1\bfw}=\overline{\Delta\bfw}=0\, ,\ \dots
$$
or, equivalently,
$$\overline{\partial_1\bfu}=\partial_1\bar \bfu,\ \overline{\Delta\bfu}=\Delta\bar\bfu\, ,\ \dots
$$

\par
For a solution of \eqref{per0} in the class $\Gamma$, we have
\be
\begin{split}
&\zeta\partial_\tau\bfw -\Delta\bar{\bfu}-\Delta\bfw+\nabla \bar{p}+\nabla{\sfp}
\\&=\lambda\,[\partial_1\bar{\bfu} +\partial_1\bfw-\bfu_0\cdot(\nabla\bar{\bfu}+\nabla\bfw)
+(\zeta\dot{\bfxi}-\bar{\bfu}-\bfw)\cdot(\nabla\bfu_0+ \nabla\bar{\bfu}+\nabla\bfw)]
 \end{split}
\eeq{perclassC0}
and
\be
\zeta^2\ddot{\bfxi}+\omega_{\sf n}^2(\bar\bfchi+\bfxi)+\varpi\Int{\partial\Omega}{} \mathbb T(\bar{\bfu},\bar{p})\cdot\bfn+\varpi\Int{\partial\Omega}{} \mathbb T(\bfw,{\sfp})\cdot\bfn=\0\,, \ \ \text{ in } (0,2\pi).
\eeq{perclassC1}
Since we assume $\bar{\sfp}=0$, we have
$$\int_0^{2\pi}\Int{\partial\Omega}{} \mathbb {\sfp} \mathbb I\cdot\bfn=\Int{\partial\Omega}{} \bar{\sfp} \mathbb I\cdot\bfn=\bf0.$$
The facts that $\bar\bfxi=\bf0$ and $\bar\bfw=\bf0$ then imply that \eqref{perclassC1} splits in
\bestar
\begin{split}
 \omega_{\sf n}^2\bar\bfchi+\varpi\Int{\partial\Omega}{} \mathbb T(\bar{\bfu},\bar{p})\cdot\bfn=\0 ,\\
\zeta^2\ddot{\bfxi}+\omega_{\sf n}^2\bfxi+\varpi\Int{\partial\Omega}{} \mathbb T(\bfw,{\sfp})\cdot\bfn =\0\,, \ \ \text{ in } (0,2\pi).
\end{split}
\eestar
Taking the mean of \eqref{perclassC0}, we infer that
$$
-\Delta\bar{\bfu}+\nabla\bar{p}=\lambda\Big[\partial_1\bar{\bfu}-\bfu_0\cdot\nabla\bar{\bfu}-\bar{\bfu}\cdot(\nabla\bfu_0+\nabla\bar{\bfu})
+ \bar{(\zeta\dot{\bfxi}-\bfw)\cdot\nabla\bfw}\, \Big].
$$
Summing up, and setting
$$\bfM(x):=\bar{(\zeta\dot{\bfxi}-\bfw)\cdot\nabla\bfw},$$
we infer that \eqref{per0} can be split into the coupled system
\be\ba{cr}\medskip\left.\ba{lr}\medskip
&-\Delta\bar{\bfu}+\nabla\bar{p}=\lambda[\partial_1\bar{\bfu}-\bfu_0\cdot\nabla\bar{\bfu}-\bar{\bfu}\cdot\nabla\bfu_0-\bar{\bfu}\cdot\nabla\bar{\bfu}
+\bfM(x)]\\ \medskip
& \Div\bar{\bfu}=0\ea\right\}\ \ \mbox{in $\Omega$\,,}\\ \medskip
\bar{\bfu}=\0\ \ \text{on }\partial\Omega\,,\\
\omega_{\sf n}^2\bar{\bfchi}+\varpi\Int{\partial\Omega}{}\mathbb T(\bar{\bfu},\bar{p})\cdot\bfn=\0\,,
\ea
\eeq{per1}
and
\be\ba{cr}\medskip\left.\ba{lr}\medskip
&\zeta\partial_\tau\bfw-\Delta\bfw+\nabla{\sfp}
 =\lambda\Big[\partial_1\bfw-\bfu_0\cdot\nabla\bfw+(\zeta\dot{\bfxi}-\bfw)\cdot\nabla\bfu_0\Big]\\ \medskip
& +\lambda[-\bfM(x)+{(\zeta\dot{\bfxi}-\bfw)\cdot\nabla\bfw}-\bar{\bfu}\cdot\nabla\bfw+(\zeta\dot{\bfxi}-\bfw)\cdot\nabla\bar{\bfu}]\\
 & \Div\bfw=0\ea\right\}\ \ \mbox{in $\Omega\times [0,2\pi]$}\\ \medskip
\bfw=\zeta\dot{\bfxi}\ \ \mbox{on }\partial\Omega\times[0,2\pi],\\
{\zeta^2}\ddot{\bfxi}+\omega_{\sf n}^2\bfxi+\varpi\Int{\partial\Omega}{}\mathbb T(\bfw,{\sfp})\cdot\bfn=\0\,,\ \ \mbox{in $[0,2\pi]$}\,.
\ea
\eeq{per2}

In  the proof of the main finding of this subsection given below, we need to use a specific ``cut-off'' function, whose properties are collected in the following lemma.
\Bl There exists $\psi_R\in C_0^\infty(\real^3)$, $R\in(0,\infty)$,  with the following properties
\begin{itemize}
  \item [(i)] $\psi_R(x)\in [0,1]$, for all $x\in\real^3$ and $R>0$;
  \item [(ii)] $\Lim{R\to\infty}\psi_R(x)=1$ for all $x\in\real^3$;
  \item [(iii)] $\psi_R(x)=1$ for all $x\in B_R$, and $\supp(\psi_R)\subset B_{2R^2}$, $R\ge 1$;
  \item [(iv)]  $\supp(\nabla\psi_R)\subset B_{2R^2}\backslash B_R=:S_R$, $R\ge 1$;
 \item [(v)] $\|\bfu|\nabla\psi_R|\|_2\le c\,\|\nabla\bfu\|_{2,\Omega^{\frac{R}{\sqrt2}}}$,\, with $c$ independent of $R$;
  \item [(vi)] $\partial_1\psi_R\in L^{2}(\Omega)$\,.
\end{itemize}
\EL{cutoff}
{\it Proof.} Let  $\psi=\psi(r)$, $r\in (0,\infty)$, be a smooth, non-increasing function that is 1 if $r<1/2$ and 0 if $r>1$, and set
$$
\psi_R(x):= \psi\left(\frac12 \sqrt{\frac{x_1^2}{R^4}+\frac{\rho^2}{R^2}}\right)\,,\ \ \rho^2:=x_2^2+x_3^2\,.
$$
We thus have
\be
\psi_R(x)=\left\{\ba{ll}\medskip 1&\, \mbox{if\, $\frac{x_1^2}{R^4}+\frac{\rho^2}{R^2}$}\le1\\
0&\, \mbox{if\, $\frac{x_1^2}{R^4}+\frac{\rho^2}{R^2}\ge4$}\ea\right.\,,
\eeq{Pgc}
which at once implies the validity of properties (i) and (ii). Furthermore, for $R\ge1$, we get
$$
\frac1{R^4}\left( x_1^2+\rho^2\right)\le \frac{x_1^2}{R^4}+\Frac{\rho^2}{R^2}\le\frac1{R^2}\left( x_1^2+\rho^2\right)\,.
$$
The latter, combined with \eqref{Pgc}, proves the statements in (iii) and (iv). Finally, the remaining properties (v) and (vi) are obtained  exactly like in \cite[Lemma II.6.4]{Gab}.\hfill$\square$\par\smallskip
The following result holds.
\Bt
Let $(\bfu,p,\bfchi)\in\Gamma$ be a solution of \eqref{per0}. If $\lambda<\lambda_2$, necessarily $(\bfu,p,\bfchi)\equiv(\0,0,\0)$.
\ET{noper}
{\it Proof.}
Recall that we use the decomposition
$$\bfu=\bar{\bfu}+\bfw,\quad p =\bar{p}+{\sfp},\quad \bfchi=\bar{\bfchi}+\bfxi.$$
We test both sides of \eqref{per2}$_1$ by $\bfw$, and integrate by parts over $\Omega\times(0,2\pi)$. Taking into account the summability properties of elements in the class $\Gamma$ and the definition \eqref{sfaco2} of $\lambda_2$, we obtain
\be
\bar{\|\nabla\bfw\|_2^2} =\lambda\left(\bar{(\zeta\dot{\bfxi}-\bfw)\cdot\nabla(\bfu_0+\bar{\bfu}),\bfw}\right)
\le \frac{\lambda}{\lambda_2}\bar{\|\nabla\bfw\|_2^2}-\lambda\left(\bfM(x),\bar{\bfu}\right)\,.
\eeq{per3}
Testing both sides of \eqref{per1}$_1$ by $\psi_R\bar{\bfu}$, with $\psi_R$ given in \lemmref{cutoff}, and integrating by parts over $\Omega$ as needed, we get
\be\ba{ll}\medskip
\|\psi_R^{\frac12}\nabla\bar{\bfu}\|_2^2=\half\lambda\,[-(\partial_1\psi_R\bar{\bfu},\bar{\bfu})+(|\bar{\bfu}|^2(\bfu_0+\bar{\bfu}),\nabla\psi_R)+2(\psi_R\bfM(x),\bar{\bfu})]
\\
\hspace*{3cm}
-\lambda(\psi_R\bar{\bfu}\cdot\nabla\bfu_0,\bar{\bfu})+(\bar{p}\,\bar{\bfu},\nabla\psi_R):=\Sum{k=1}5\,I_k\,.
\ea
\eeq{per4}
Using (iv)--(vi) in \lemmref{cutoff}  along with H\"older inequality, we show
\bestar
|I_1|+|I_2| \le \half\lambda\,\left(\|\partial_1\psi_R\|_2+\|(\bfu_0+\bar{\bfu})|\nabla\psi_R|\|_2\right)\|\bar{\bfu}\|_{4,S_R}^2
\le c\,(1+\|\nabla\bfu_0\|_2+\|\nabla\bar{\bfu}\|_2)\|\bar{\bfu}\|_{4,S_R}^2
\eestar
and also
$$
|I_5|\le \|\bar{\bfu}|\nabla\psi_R|\|_2\|\bar{p}\|_{2,S_R}\le c\,\|\nabla\bar{\bfu}\|_2\|\bar{p}\|_{2,S_R}\,.
$$
We now pass to the limit $R\to\infty$ in \eqref{per4}. With the help of the last two displayed inequalities along with \lemmref{cutoff}-(ii),  the fact that $\bfu,\bfw\in \Gamma$ and the definition of $\lambda_1$, we obtain
\be
\|\nabla\bar{\bfu}\|_2^2 =-\lambda\,(\bar{\bfu}\cdot\nabla\bfu_0,\bar{\bfu})+\lambda(\bfM(x),\bar{\bfu})
\le\displaystyle\frac\lambda{\lambda_1}\|\nabla\bar{\bfu}\|_2^2+\lambda(\bfM(x),\bar{\bfu})\,.
\eeq{per5}
Summing side-by-side \eqref{per3} and \eqref{per5}, we infer that
$$
\left(1-\frac{\lambda}{\lambda_2}\right)\bar{\|\nabla\bfw\|_2^2}
+\left(1-\frac{\lambda}{\lambda_1}\right) \|\nabla\bar{\bfu}\|_2^2
\le 0,
$$
 and since $\frac\lambda{\lambda_2}<1$ implies $\frac\lambda{\lambda_1}<1$, we conclude that
 $$
\nabla{\bfu}(x,t)= \0\,,\ \ \mbox{for all $(x,t)\in \Omega\times [0,2\pi]$\,,}
$$
which, by \eqref{1.9} and \eqref{1.10} concludes the proof
of the theorem.
\hfill$\square$\par\smallskip
\setcounter{equation}{0}
\subsection{Steady-State Bifurcation}\label{steady}
We shall now undertake the study of loss of uniqueness and occurrence of bifurcation for the family of solutions ${\sf s}_0(\lambda)$, $\lambda>0$, whose existence was established in \theoref{exi}. More precisely, in this subsection we will furnish necessary and sufficient conditions for the occurrence of steady bifurcation.
\par
Let $\lambda_s>0$, $U(\lambda_s)$ be a neighborhood of $\lambda_s$, and denote by
\be
{\sf s}_0(\lambda):=(\bfu_0(\lambda),p_0(\lambda),\bfchi_0(\lambda))\,,\ \ \lambda\in U(\lambda_s)\,,
\eeq{stbr}
a first solution to \eqref{03} determined in \theoref{exi}. Next, let
$$(\bfu_0(\lambda)+\bfu(\lambda),{p}_0(\lambda)+p(\lambda),\bfchi_0(\lambda)+\bfchi(\lambda))\,,\ \ \lambda\in U(\lambda_s)\,,
$$
be another solution to \eqref{03} so that $(\bfu(\lambda),{p}(\lambda),\bfchi(\lambda))$ solves the following homogeneous boundary-value problem
\be\ba{cc}\medskip\left.\ba{lr}\medskip
&-\Delta\bfu+\nabla p=\lambda\,(\partial_1\bfu-\bfu_0\cdot\nabla\bfu-\bfu\cdot\nabla\bfu_0-\bfu\cdot\nabla\bfu)\\
&\Div\bfu=0\ea\right\}\ \ \mbox{in $\Omega$}\,,\\ \medskip
\bfu(x)=\0\,, \ \mbox{ $x\in\partial\Omega$}\,;\ \
\Lim{|x|\to\infty}\bfu(x)=\0\,,\\
\omega_{\sf n}^2\bfchi = -\varpi\Int{\partial\Omega}{} \mathbb T(\bfu,p)\cdot\bfn\,.
\ea
\eeq{pErt}
Then, formally, steady-state bifurcation reduces to show that\begin{itemize} \item[(i)] $(\bfu(\lambda),{p}(\lambda),\bfchi(\lambda))\not\equiv (\0,0,\0)$\,,\ \ $\lambda\in U(\lambda_s)$\,,\item[(ii)] $(\bfu(\lambda),{p}(\lambda),\bfchi(\lambda))\to (\0,0,\0)$ as $\lambda\to\lambda_s$\,,\end{itemize}
in which case, $\left(\lambda_s,{\sf s}_0(\lambda_s)\right)$ is called a bifurcation point for problem \eqref{03} (or, equivalently, $\left(\lambda_s,\0\right)$ a bifurcation point for \eqref{pErt}).
The above properties should be, of course, rigorously formulated  and their validity properly ascertained in an  appropriate functional setting.
\Br We observe that, in order to prove (i)-(ii) it is enough  to prove
\begin{itemize} \item[(i)$^\prime$] $(\bfu(\lambda),{p}(\lambda))\not\equiv (\0,0)$\,,\ \ $\lambda\in U(\lambda_s)$\,,\item[(ii)$^\prime$] $(\bfu(\lambda),{p}(\lambda))\to (\0,0)$ as $\lambda\to\lambda_s$\,,\end{itemize}
{\it provided} the existence of the bifurcating branch is obtained in a function class that allows the control of $\bfchi(\lambda)$, namely, the right-hand side of \eqref{pErt}$_5$. We could thus restrict ourselves to prove such a type of existence. However, also in view of the analysis of time-periodic bifurcation that we shall develop in \cite{BoGaGaper}, we prefer to study problem \eqref{pErt} as a whole.
\ER{rem}

Let
\be
\calz:=[D^{2,2}(\Omega)\cap\cald_0^{1,2}(\Omega)]\times\real^3\,,
\eeq{z}
and
\be
\mathcal X:= X^{2}(\Omega)\times\real^3\,;
\eeq{x}
see \eqref{ics}.
Clearly,
$\mathcal X\subset\calz$.
Given a function $\bfu:\Omega\to\real^3$ and a vector $\bfchi\in\real^3$, we put
$$
\bsfU:=(\bfu,\bfchi)\,.
$$
Our first goal to analyse steady bifurcation is to rewrite the left-hand side of \eqref{pErt}$_1$ with suitable operators acting on $\bsfU$ (that will therefore also include the compatibility condition \eqref{pErt}$_5$ for $\bfchi$).
To this aim, we consider several maps and their properties. Define first
$$
\hat{\Delta}:\bsfU\in \calz\mapsto \hat{\Delta}(\bsfU)\in Y\,,
$$
see \eqref{ips}, where
\be
\hat{\Delta}(\bsfU)=\left\{\ba{ll}\medskip
 -\Delta\bfu\ \mbox{in $\Omega$}\,,\\
\omega_{\sf n}^2\bfchi+2\varpi\Int{\partial\Omega}{}\mathbb D(\bfu)\cdot\bfn\ \ \mbox{in $\Omega_0$\,,}\ea\right.
\eeq{Deltah}
and set $\hat{\mathscr A}:=\mathscr P\,\hat{\Delta}$, where $\mathscr P$ is the self-adjoint orthogonal projection of $\call^2(\real^3)$ onto $\calh(\real^3)$.  We thus have
$$
\hat{\mathscr A}:\bsfU\in \calz \mapsto \hat{\mathscr A}(\bsfU)\in \caly\,,
$$
see \eqref{ips}, and, by \eqref{Deltah} and \lemmref{0.1},
\be
\hat{\mathscr A}(\bsfU)=\left\{\ba{ll}\medskip
 -\Delta\bfu+\nabla p\ \mbox{in $\Omega$}\,,\\
\omega_{\sf n}^2\bfchi+\varpi\Int{\partial\Omega}{}\left(2\mathbb D(\bfu)-p\,\mathbb I\right)\cdot\bfn\ \ \mbox{in $\Omega_0$\,,}\ea\right.,
\eeq{Stoke}
for some $p\in\mathcal G(\Omega)$.\par
The following result holds.
\Bl The operator $\hat{\mathscr A}:\calz\to\caly$ is a homeomorphism.
\EL{Sto}
{\it Proof.} We have to show that, for any $(\bff,\bfF)\in [\cald_0^{-1,2}(\Omega)\cap L^2(\Omega)]\times\real^3$ there exists one and only one $(\bfu,\bfchi)\in [D^{2,2}(\Omega)\cap\cald_0^{1,2}(\Omega)]\times\real^3$ satisfying
\be\ba{ll}\medskip
-\Delta\bfu+\nabla p=\bff\,,\ \ \Div\bfu=0\ \mbox{in $\Omega$}\,;\ \ \bfu=\0\ \ \mbox{on $\partial\Omega$}\,,\\
\omega_{\sf n}^2\bfchi+\varpi\Int{\partial\Omega}{}\left(2\mathbb D(\bfu)-p\,\mathbb I\right)\cdot\bfn=\bfF\,.\ea
\eeq{Sto1}
The result will then follow by the open mapping theorem.
It is well-known that for any $\bff\in~ \cald_0^{-1,2}(\Omega)$ there exists a unique $(\bfu,p)\in \cald_0^{1,2}(\Omega)\times L^2(\Omega)$ satisfying \eqref{Sto1}$_{1,2,3}$ in distributional sense, see for instance \cite[Theorem V.2.1]{Gab}. Moreover, since $\bff\in L^2(\Omega)$, we also have $(\bfu,p)\in D^{2,2}(\Omega)\times D^{1,2}(\Omega)$ \cite[Theorems IV.5.1 and V.5.3]{Gab}. By the trace theorem, for a fixed $R>R_*$,
$$
\Big|\Int{\partial\Omega}{}\left(2\mathbb D(\bfu)-p\,\mathbb I\right)\cdot\bfn\Big|\le c\,(\|\bfu\|_{2,2,\Omega_{R}}+\|p\|_{1,2,\Omega_{R}})<\infty\,,
$$
so that $\bfchi$ is uniquely determined from \eqref{Sto1}$_4$.
\hfill$\square$\par\smallskip
Let
$$
\hat{\partial_1}:\bsfU\in \calx\mapsto \hat{\partial_1}(\bsfU)\in \calh(\real^3)
$$
with
\be
\hat{\partial_1}(\bsfU)=\left\{\ba{ll}\medskip
 -\partial_1\bfu\ \mbox{in $\Omega$}\,,\\
\0\ \ \mbox{in $\Omega_0$\,.}\ea\right.
\eeq{D1h}
It is readily checked that, as stated, $\hat{\partial}_1(\bsfU)\in \calh(\real^3)$. By \lemmref{0.1} and \eqref{D1h}, this amounts to show that
\be
\int_{\Omega}\partial_1\bfu\cdot\nabla p=0\,,\ \ \mbox{for all $p\in D^{1,2}(\Omega)$}\,.
\eeq{wwh}
Since $\bfu\in\cald^{1,2}_0(\Omega)$, there is a sequence $(\bfu_k)_k\subset \calc_0(\Omega)$ such that $\|\nabla(\bfu_k-\bfu)\|_2\to 0$ as $k\to\infty$. Then, by an integration by parts combined with the condition $\Div\bfu_k=0$, we deduce that \eqref{wwh} holds for each $\bfu_k$ and hence for $\bfu$ after passing to the limit $k\to\infty$.
\par
Next, let
$$
\mathscr C:\bsfU\in \mathcal X\mapsto \mathscr C(\bsfU)\in Y
$$
where
\be
\mathscr C(\bsfU)=\left\{\ba{ll}\medskip
 \bfu_0\cdot\nabla\bfu+\bfu\cdot\nabla\bfu_0\ \mbox{in $\Omega$}\,,\\
\0\ \ \mbox{in $\Omega_0$\,,}\ea\right.
\eeq{Kh}
and consider the following operator
\be
\mathscr L_\lambda:\bsfU\in \mathcal X\mapsto\hat{\mathscr A}(\bsfU) +\lambda[\hat{\partial_1}(\bsfU)+\mathscr P\mathscr C(\bsfU)]\in \caly\,.
\eeq{2.2lin}

Then we prove

\Bl  For any $\lambda\neq 0$, the operator $\mathscr L_\lambda:\mathcal X\to\caly$ is Fredholm of index 0.
\EL{fredo}
{\it Proof.}
We first show that $\mathscr L_\lambda$ is well defined.
Clearly $\mathscr P\bfphi=\bfphi$,  for  $\bfphi\in \mathcal C_0(\Omega)$. Therefore, for any such a $\bfphi$, integrating by parts we get
$$
\left(\mathscr L_\lambda(\bsfU),\bfphi\right)=\lambda(\partial_1\bfu,\bfphi)-\left(\lambda(\bfu_0\otimes\bfu+\bfu\otimes\bfu_0)+\nabla\bfu,\nabla\bfphi\right).
$$
Since $\bfu\in X(\Omega)$, employing  H\"older inequality and \lemmref{3.1} we deduce that
$\mathscr L_\lambda(\bsfU)\in\cald_0^{-1,2}(\Omega)$.
Moreover,
$$
\|\bfu_0\cdot\nabla\bfu+\bfu\cdot\nabla\bfu_0\|_2\le \|\bfu_0\|_\infty\|\nabla\bfu\|_2+\|\bfu\|_\infty\|\nabla\bfu_0\|_2\,.
$$
Since $\bfu\in X^2(\Omega)$,
 by \lemmref{Xemb} and \theoref{exi}, we conclude that $\mathscr L_\lambda(\bsfU)\in \calh(\real^3)$ as claimed.

We now turn to the verification of the Fredholmness of $\mathscr L_\lambda$.
We decompose $\mathscr L_\lambda$ as follows
$$\mathscr L_\lambda=\mathscr L_\lambda^0+\lambda\mathscr P\mathscr C,$$
where $\mathscr L_\lambda^0:=\hat{\mathscr A}+\lambda\,\hat{\partial_1}.$
To prove the Fredholm property for $\mathscr L_\lambda$, it is enough to show that
$\mathscr L_\lambda^{0}$ is a homeomorphism and that $\mathscr C$ is a compact operator.

We start by showing $\mathscr L_\lambda^{0}$ is a homeomorphism, that is, for any $(\bff,\bfF)\in [\cald_0^{-1,2}(\Omega)\cap L^2(\Omega)]\times\real^3$ there exists unique $(\bfu,\bfchi)\in X^2(\Omega)\times \real^3$ solving
\be\ba{cc}\medskip\left.\ba{ll}\medskip
\lambda\,\partial_1\bfu+\Delta\bfu=\nabla p+\bff\\ \medskip
\Div\bfu=0\ea\right\}\ \ \mbox{in $\Omega$}\,,\\ \medskip
\bfu(x)=\0\ \ \mbox{on $\partial\Omega$}\,,\\
\omega_{\sf n}^2\bfchi+\varpi\Int{\partial\Omega}{}\mathbb T(\bfu,p)\cdot\bfn=\bfF\,.
\ea
\eeq{chia_p}
To show the validity of this property,
we notice that corresponding to the given $\bff$, from \cite[Theorem 2.1]{GaEi} we know that there exists a unique solution $(\bfu,p)\in X^2(\Omega)\times W^{1,2}(\Omega)$ to \eqref{chia_p}$_{1,2,3}$.
By the trace theorem we get
\be
\left|\int_{\partial\Omega}\mathbb T(\bfu,p)\cdot\bfn\right|\le c\,(\|\bfu\|_{2,2,\Omega_R}+\|p\|_{1,2,\Omega_R})<\infty\,.
\eeq{2.3_s}
The associated and uniquely determined displacement $\bfchi$ is then obtained from \eqref{chia_p}$_4$,
which concludes the proof of the homeomorphism property of $\mathscr L_\lambda^0$ by the open mapping theorem.

We next show that the operator $\mathscr C$ is compact. To this end, consider a sequence
$$\{\bsfU_k\}_k\equiv \{(\bfu_k,\bfchi_k)\}\subset \calx$$
with
\be
\|\bfu_k\|_{X^2}+|\bfchi_k|\le M\,,
\eeq{unibou}
and $M$ independent of $k\in\nat$.
This implies the existence of $\bfu_*\in L^4(\Omega)\cap D^{1,2}(\Omega)\cap D^{2,2}(\Omega)$ and  $\bfchi_*\in\real^3$ such that along a subsequence (that we continue to denote by $\{(\bfu_k,\bfchi_k)\}$)
\be\ba{ll}\medskip
\bfchi_k-\bfchi_*\to 0\,,\ \ \mbox{in $\real^3$}\\ \medskip
\bfw_k:=\bfu_k-\bfu_*\to 0\,,\ \ \mbox{weakly in $L^4(\Omega)$}\,,\\
\nabla\bfw_k\to 0 \,,\ \ \mbox{weakly in $W^{1,2}(\Omega)$}\,.
\ea
\eeq{weakc}
Moreover, by \lemmref{Xemb} and compact embedding results, we also have
\be
\bfw_k\to 0\,,\ \ \mbox{strongly in $W^{1,2}(\Omega_R)$, for all $R>R_*$.}
\eeq{strongc}
Observe that
$$
|\bfu_0\cdot\nabla\bfw_k+\bfw_k\cdot\nabla\bfu_0|_{-1,2}\le\|\bfu_0\otimes\bfw_k\|_2\le  \|\bfu_0\|_\infty\|\bfw_k\|_{2,\Omega_R}+\|\bfu_0\|_{4,\Omega^R}\|\bfw_k\|_4\,,
$$
where we recall $\Omega^R=\Omega\backslash\bar{\Omega_R}=\real^3\setminus\bar{B_R}$,
and
$$
\|\bfu_0\cdot\nabla\bfw_k+\bfw_k\cdot\nabla\bfu_0\|_{2}\le \|\bfu_0\|_{1,\infty}\|\bfw_k\|_{1,2,\Omega_R}+\|\bfu_0\|_{1,4,\Omega^R}\|\bfw_k\|_{1,4}\,.
$$
From \theoref{exi}, we know that $\bfu_0\in W^{1,\infty}(\Omega)\cap W^{1,4}(\Omega)$, so that letting first $k\to\infty$ in the above two inequalities and using \eqref{unibou}--\eqref{strongc}, and then $R\to\infty$, we deduce
$$
\lim_{k\to\infty}\left(|\bfu_0\cdot\nabla\bfw_k+\bfw_k\cdot\nabla\bfu_0|_{-1,2}+\|\bfu_0\cdot\nabla\bfw_k+\bfw_k\cdot\nabla\bfu_0\|_{2}\right)=0\,.
$$
In view of \eqref{Kh}, this proves that $\mathscr C$ is compact and finishes the proof. \hfill$\square$\par\smallskip

From \theoref{exi}, we know that the  solution ${\sf s}_0(\lambda):=(\bfu_0(\lambda),p_0(\lambda),\bfchi_0(\lambda))$ to
problem \eqref{03} is unique, provided that $\lambda<\lambda_1(\lambda)$, as defined in \eqref{sfaco1}. Our next concern
is to furnish a sufficient condition for  local  uniqueness of ${\sf s}_0(\bar{\lambda})$ whenever $\bar{\lambda}\ge \lambda_1(\bar{\lambda})$ if this situation occurs. Let
$$(\bsfu_0,{\sfp}_0,{\sf \bftau}_0):=(\bfu_0(\bar{\lambda}),p_0(\bar{\lambda}),\bfchi_0(\bar{\lambda}))$$
be a given solution of the problem
\be\ba{cc}\medskip\left.\ba{ll}\medskip
-\Delta\bfu_0+\nabla p_0=\bar\lambda\,(\partial_1\bfu_0-\bfu_0\cdot\nabla\bfu_0)\\
\Div\bfu_0=0\ea\right\}\ \ \mbox{in $\Omega$}\,,\\ \medskip
\bfu_0(x)=\bfe_1\,, \ \mbox{ $x\in\partial\Omega$}\,;\ \
\Lim{|x|\to\infty}\bfu_0(x)=\0\,,\\
\omega_{\sf n}^2\bfchi_0+\varpi\Int{\partial\Omega}{} \mathbb T(\bfu_0,p_0)\cdot\bfn=\0\,.
\ea
\eeq{03-bar}
We will write any solution $(\bfu_0(\lambda),p_0(\lambda),\bfchi_0(\lambda))$ to \eqref{03} as
$$(\bfu(\lambda)+\bsfu_0,p(\lambda)+{\sfp}_0,\bfchi(\lambda)+\bftau_0).$$
Then $(\bfu,p,\bfchi)=(\bfu(\lambda),p(\lambda),\bfchi(\lambda))$ solves
\be\ba{cr}\medskip\left.\ba{lr}\medskip
&-\Delta\bfu+\nabla p-\bar\lambda\partial_1\bfu+\bar\lambda\left(\bsfu_0\cdot\nabla\bfu
+\bfu\cdot\nabla\bsfu_0\right)\\
&=\mu\Big(\partial_1(\bsfu_0+\bfu)-(\bsfu_0+\bfu)\cdot\nabla(\bsfu_0+\bfu)-\bar\lambda\mu^{-1} \bfu\cdot\nabla\bfu\Big)\\
& \Div\bfu=0\ea\right\}\ \ \mbox{in $\Omega$}\,,\\
\medskip
\bfu(x)=\0\,, \ \mbox{ for }x\in\partial\Omega\,, \Lim{|x|\to\infty}\bfu(x)=\0 \,,\\
\omega_{\sf n}^2\bfchi+\varpi\Int{\partial\Omega}{} \mathbb T(\bfu,p)\cdot\bfn=\0\,,
\ea
\eeq{03_1111}
where we have set $\mu:=\lambda-\bar\lambda$. To rewrite \eqref{03_1111} in a suitable way, we define one more operator:
$$
\mathscr O:\bsfU\in \mathcal X\mapsto \mathscr O(\bsfU)\in Y
$$
where
\be
{\mathscr O}(\bsfU)=\left\{\ba{ll}\medskip
 -\partial_1(\bsfu_0+\bfu)  + (\bsfu_0+\bfu)\cdot\nabla(\bsfu_0+\bfu)+\bar{\lambda}\,\mu^{-1}\bfu\cdot\nabla\bfu\ \ \mbox{in $\Omega$}\,,\\
\0\ \ \mbox{in $\Omega_0$\,.}\ea\right.
\eeq{M}
It is clear that \eqref{03_1111} is formally equivalent to
$$\mathscr L_{\bar\lambda}(\bsfU)+\mu\mathscr P\mathscr O(\bsfU)=0.$$
We claim that the map $\mathscr O$ is well defined. Indeed, taking into account \theoref{exi} and the fact that $\bfu\in X^2(\Omega)$, by arguments similar to those employed previously we easily show that, with $\bfa,\bfb$ being either $\bsfu_0$ or $\bfu$,
\be
\partial_1\bfu, \,  \bfa\cdot\nabla\bfb\in  \cald_0^{-1,2}(\Omega)\cap L^2(\Omega)\,,\ \  \ \partial_1\bsfu_0\in L^2(\Omega)\,.
\eeq{mcsnf}
Moreover, testing \eqref{03}$_1$ (with $\bfu_0\equiv \bsfu_0$, $p_0\equiv {\sfp}_0$, $\lambda\equiv\bar \lambda$) by $\bfphi\in \calc_0(\Omega)$ and integrating by parts, we get
$$
\lambda_s(\partial_1\bsfu_0,\bfphi)=-(\bsfu_0\otimes\bsfu_0,\nabla\bfphi)+(\nabla\bsfu_0,\nabla\bfphi)\,.
$$
Since $\bsfu_0\in L^4(\Omega)\cap D^{1,2}(\Omega)$, from the latter we deduce $\partial_1\bsfu_0\in \cald_0^{-1,2}(\Omega)$, which, along with \eqref{mcsnf} proves that $\mathscr O$ is well defined.

We are now in position to give a sufficient condition for local uniqueness.

\Bt Let $(\bsfu_0,{\sfp}_0,{\sf \bftau}_0):=(\bfu_0(\bar{\lambda}),p_0(\bar{\lambda}),\bfchi_0(\bar{\lambda}))$,  and  assume that the equation
\be
\mathscr L_{\bar{\lambda}}(\bsfU)=\0
\eeq{linun}
has only the solution $\bsfU\equiv\0$. Then, there exists a neighborhood $U(\bar{\lambda})\supset\{\bar{\lambda}\}$ such that for $\lambda\in U(\bar{\lambda})$, ${\sf s}_0(\lambda):=(\bfu(\lambda),p_0(\lambda),\bfchi_0(\lambda))$ is the only solution to \eqref{03}. Moreover, $\lambda\mapsto{\sf s}_0(\lambda)$ is analytic at $\lambda=\bar{\lambda}$, and
$$
(\bfu_0(\lambda),p_0(\lambda),\bfchi_0(\lambda))\to (\bsfu_0,{\sfp}_0,{\bftau}_0)\ \ \mbox{as $\lambda\to\bar{\lambda}$\,.}
$$
\ET{2.1}
{\it Proof.}
Let
$$
\mathcal F:(\bsfU,\mu)\in  \Xi(0) \times \cali(0)\mapsto \mathscr L_{\bar \lambda}(\bsfU)+\mu\mathscr P\mathscr O(\bsfU)\in\mathcal Y\,,
$$
where $\Xi(0)\times\cali(0)$ is a neighborhood of $({\bf0},0)\in \mathcal X\times \real$.
We have seen that \eqref{03} translates in \eqref{03_1111} which is equivalent to
\be
\mathcal F(\bsfU,\mu)=\0\in \mathcal Y\,.
\eeq{2.5}
Clearly, \eqref{2.5} has the solution $(\bsfU=\0, \mu=0)$. Moreover, it is standard to verify that $\mathcal F$ is Frech\'et-differentiable with derivative $D_{{\mbox{\tiny $\bsfU$}}}\mathcal F(\0,0)=\mathscr L_{\bar{\lambda}}$. Since $\mathscr L_{\bar{\lambda}}$ is Fredholm of index 0, the assumption made in \eqref{linun}  implies that $\mathscr L_{\bar{\lambda}}$ is a homeomorphism. In addition, $\calf$ is polynomial in $(\bsfU,\mu)$.
As a consequence, by the analytic version of the implicit function theorem we show the  property stated in the theorem, which is thus completely established.\hfill$\square$\par\smallskip

\theoref{2.1} tells us, in particular, that ${\sf s}_0(\bar{\lambda})$, $\bar{\lambda}\ge \lambda_1$, is unique as long as the corresponding linearization satisfies  \eqref{linun}. Moreover, this solution can be analytically (and uniquely) continued up to the first value of $\bar{\lambda}$, say, $\lambda_s$, where  \eqref{linun} is violated.\par
We thus have:
\Bc A necessary condition for $(\lambda_s,{\sf s}_0(\lambda_s))$ to be a steady-state bifurcation point is
\be
{\rm dim}\,{\sf N}[\mathscr L_{\lambda_s}]>0\,.
\eeq{NoNuni}
\EC{6.1}
\Br Recalling the definition of $\mathscr L_\lambda$ in \eqref{2.2lin}, we have that condition \eqref{NoNuni} is {\it equivalent} to the request that the linear problem
\be\ba{cc}\medskip\left.\ba{ll}\medskip
-\Delta\bfu+\nabla p=\lambda_s\,(\partial_1\bfu-\bfu_0(\lambda_s)\cdot\nabla\bfu-\bfu\cdot\nabla\bfu_0(\lambda_s))\\
\Div\bfu=0\ea\right\}\ \ \mbox{in $\Omega$}\,,\\ \medskip
\bfu(x)=\0\,, \ \mbox{ $x\in\partial\Omega$}\,;\qquad\Lim{|x|\to\infty}\bfu(x)=\0\,,\\
\omega_{\sf n}^2\bfchi+\varpi\Int{\partial\Omega}{} \mathbb T(\bfu,p)\cdot\bfn=\0\,,
\ea
\eeq{LiPr}
has at least one solution $(\bfu,\bfchi)\in X^2(\Omega)\times \real\backslash \{\0,\0\}$. This is, in turn, {\it equivalent} to the condition that \eqref{LiPr}$_{1,2,3,4}$ has a non-identically zero solution $\bfu\in X^2(\Omega)$. We may then conclude that the necessary condition for bifurcation in {\it absence} of compatibility condition \eqref{LiPr}$_{5}$ remains such also in its presence; see also \remref{rem}.
\ER{6.}

Our next goal is to find sufficient conditions for the occurrence of steady bifurcation. For this, we shall rewrite \eqref{pErt} as an equation in a suitable Banach space, and establish some basic properties of the involved operators. After that, we will be able to apply abstract  bifurcation results to our case and derive the desired conditions.\par
To reach these purposes, we need to introduce another operator, namely
$$
\mathscr N:\bsfU\in \mathcal X\mapsto \mathscr N(\bsfU)\in Y
$$
where
\be
\mathscr N(\bsfU)=\left\{\ba{ll}\medskip
 \bfu\cdot\nabla\bfu\ \mbox{in $\Omega$}\,,\\
\0\ \ \mbox{in $\Omega_0$\,.}\ea\right.
\eeq{NoL}
\par
By arguing as we did for the map $\mathscr O$ in \theoref{2.1}, we show that $\mathscr N$ is well defined. Thus, from \eqref{NoL} and \eqref{2.2lin}, we deduce that problem \eqref{pErt} is equivalent to the following equation
\be
\hat{\mathscr A}(\bsfU) +\lambda[\hat{\partial_1}(\bsfU)+\mathscr P\mathscr C(\lambda)(\bsfU)+
\mathscr P\mathscr N(\bsfU)]=\0\ \ \,\mbox{in $\mathcal Y$}\,,\ \ \lambda\in U(\lambda_s),
\eeq{BiP}
where we have emphasized that the operator $\mathscr C$ depends on $\lambda$ through $\bfu_0$.
In view of \lemmref{Sto}, we may operate on both sides of \eqref{BiP} with $\hat{\mathscr A}^{-1}$, so that \eqref{BiP} becomes:
\be
{\mathscr F}(\bsfU,\lambda)=\0\,,\ \ \,\mbox{in $\mathcal Z$}\,,\ \ \lambda\in U(\lambda_s)
\eeq{BiP1}
where
\be
{\mathscr F}:={\sf I}+\lambda\mathscr M(\lambda)+\lambda \mathscr R
\eeq{BiP2}
with ${\sf I}$ identity in $\calz$, and
$$\ba{ll}\medskip
\mathscr M(\lambda):\bsfU\in \mathcal X\equiv {\sf D}[\mathscr M(\lambda)]\subset \calz\mapsto \hat{\mathscr A}^{-1}[\hat{\partial_1}(\bsfU)+\mathscr P\mathscr C(\lambda)(\bsfU)]\in \mathcal Z\,,\\
\mathscr R: \bsfU\in\mathcal X\equiv {\sf D}[\mathscr R]\mapsto \hat{\mathscr A}^{-1}\mathscr N(\bsfU)\in \calz\,.
\ea
$$
\Br
Since $\mathscr R$ is bilinear in $\bsfU$, for each fixed $\lambda$, the map ${\mathscr F}$ is analytic in $\bsfU$.
\ER{Ana}
\par
We shall now establish some important properties of the operator $\mathscr M(\lambda)$.
\Bl For each fixed $\lambda>0$, the operator $\mathscr M=\mathscr M(\lambda)$ is densely defined and closed.\footnote{In fact, the result holds for all $\lambda\neq0$, but this generalization is irrelevant to our aims.}
\EL{M}
{\it Proof.} Recalling the definitions given in \eqref{z} and \eqref{x}, the density property means that, given arbitrary $\bfu\in D^{2,2}(\Omega)\cap\cald_0^{1,2}(\Omega)$ and $\varepsilon>0$, there exists $\bfu_\varepsilon\in X^2(\Omega)\cap\cald_0^{1,2}(\Omega)$ such that
\be
\|\nabla(\bfu-\bfu_\varepsilon)\|_{1,2}<\varepsilon\,.
\eeq{denT}
Let $\phi_R$, $R>R_*$, be a smooth, non-increasing function of $|x|\in [0,\infty)$, such that $\phi_R(x)=1$, if $|x|\le R$, $\phi_R(x)=0$, if $|x|\ge 2R$, and
$$
|\nabla\phi_R(x)|\le c\,R^{-1}\,,\ \ |\nabla(\nabla\phi_R(x))|\le c\,R^{-2}
\,,\ \ \, \mbox{for all $x\in\real^3$}\,,
$$
with $c$ independent of $R$.
Consider the problem
$$\ba{ll}\medskip
\Div\bfv=\nabla\phi_R\cdot\bfu\,,\ \ \mbox{in $B_{R,2R}:=\{x\in \Omega: |x|\in (R,2R)\}$}\,,\ \ \bfv\in W^{2,2}_0(B_{R,2R}).
\ea
$$
It is well known that the field $\bfv$ exists and that there exists a positive constant $c_0$, independent of $R$, such that
\[
\|\nabla\bfv\|_2\le c_0\,\|\nabla\phi_R\cdot\bfu\|_{2}\,,\ \ \|D^2\bfv\|_2\le c_0\,\|\nabla(\nabla\phi_R\cdot\bfu)\|_{2}.
\]
We refer to  \cite[Theorem III.3.3 and Lemma III.3.3]{Gab}.
 Extending $\bfv$ by $\0$ outside $\Omega_R$, we deduce, in particular, $\bfv\in W^{2,2}_0(\Omega)$. Moreover, by H\"older inequality and the properties of $\phi_R$,
\be \ba{ll}\medskip
\|\nabla\bfv\|_2\le c R^{-1}\|\bfu\|_2\le c_2\|\bfu\|_{6,B_{R,2R}}\,,\\
\|D^2\bfv\|_2\le c \left(R^{-2}\|\bfu\|_{2,B_{R,2R}}+R^{-1}\|\nabla\bfu\|_2\right)\le c \,R^{-1} \left(\|\bfu\|_{6}+\|\nabla\bfu\|_2\right)\,,
\ea
\eeq{6.}
where, here and in the rest of the proof, $c$ denotes a positive constant independent of $R$. Setting
$$
\bfu_R:=\phi_R\bfu-\bfv\,,
$$
we establish at once with the help of \lemmref{0.1} that $\bfu_R\in X^{2}(\Omega)\cap \cald_0^{1,2}(\Omega)$. Also, again by the properties of $\phi_R$,  H\"older inequality and \eqref{6.} we show, in a similar manner,
\be\ba{ll}\medskip
\|\nabla(\bfu-\bfu_R)\|_2\le \|(1-\phi_R)\nabla\bfu\|_2+c \|\bfu\|_{6,B_{R,2R}}\\
\|D^2(\bfu-\bfu_R)\|_2\le \|(1-\phi_R)D^2\bfu\|_2+c R^{-1}\left(\|\bfu\|_{6}+\|\nabla\bfu\|_2\right)\,.
\ea
\eeq{6..}
Since, by \lemmref{0.1}, $\bfu\in L^6(\Omega)$, \eqref{denT} follows from \eqref{6..}, by taking $R$ sufficiently large.

To prove $\mathscr M$ is closed, we take $\{\bsfU_k\equiv(\bfu_k,\bfchi_k)\}\subset \calx$,  $\left(\bsfU\equiv(\bfu,\bfchi),\bsfV\equiv(\bfv,\bfzeta)\right)\in \calz$ such that
\be\ba{ll}\medskip
\|\nabla(\bfu_k-\bfu)\|_{1,2}+|\bfchi_k-\bfchi|\to 0\\
\mathscr M(\bsfU_k)\to \bsfV\ \ \mbox{in $\calz$}\ea\ \ \mbox{as $k\to\infty$}\,,
\eeq{CoN}
where, for simplicity, we have omitted the dependence of $\mathscr M$ on $\lambda$.
We need to show
\begin{itemize}
  \item [(a)] $\bfu\in X(\Omega)$\,;
  \item [(b)] $\mathscr M(\bsfU)=\bsfV$\,.
\end{itemize}
Let $(\bfv_k,\bfzeta_k)\equiv \bsfV_k:=\mathscr M(\bsfU_k)$, and set ${\bfv}_{kk'}=\bfv_k-\bfv_{k'}$, $\bfzeta_{kk'}=\bfzeta_k-\bfzeta_{k'}$, ${\bfu}_{kk'}=\bfu_k-\bfu_{k'}$ and  $p_{kk'}=p_k-p_{k'}$. Then,
\be\ba{cc}\medskip\left.\ba{ll}\medskip
-\Delta \bfv_{kk'}+\nabla p_{kk'}=\lambda\big(\partial_1\bfu_{kk'}-\bfu_0\cdot\nabla\bfu_{kk'}-\bfu_{kk'}\cdot\nabla\bfu_0\big)\\
\Div \bfu_{kk'}=0\ea\right\}\ \ \mbox{in $\Omega$\,,}\\ \medskip
\bfu_{kk'}=\0\ \ \mbox{on $\partial\Omega$\,,}\\
\omega_{\sf n}^2\bfzeta_{kk'}+\varpi\Int{\partial\Omega}{}\mathbb T(\bfv_{kk'},p_{kk'})\cdot\bfn=\0\,.
\ea
\eeq{ClO}
Testing \eqref{ClO}$_1$ with $\bfphi\in\calc_0(\Omega)$, integrating by parts, and using H\"older inequality we infer that
$$\ba{rl}\medskip
\lambda|(\partial_1\bfu_{kk'},\bfphi)|=&\!\!\!|\left(\nabla\bfv_{kk'}-\lambda(\bfu_0\otimes\bfu_{kk'}+\bfu_{kk'}\otimes\bfu_0),\nabla\bfphi\right)|\\
\le&\!\!\!\left(\|\nabla\bfv_{kk'}\|_2+2\lambda\|\bfu_0\|_3\|\bfu_{kk'}\|_6\right)\|\nabla\bfphi\|_2
\,.
\ea
$$

Now, $\|\nabla\bfv_{kk'}\|_2\to 0$ (by \eqref{CoN}$_2$),  $\|\bfu_{kk'}\|_6\to 0$ (by \eqref{1.9} and \eqref{CoN}$_1$) and $\bfu_0\in L^3(\Omega)$ (by \theoref{exi}), and so we conclude that $(\bfu_k)_k$ is converging also in $X(\Omega)$, which proves (a).\par
Writing \eqref{ClO} with $\bfv_{kk'}\equiv\bfv_k$, $\bfu_{kk'}\equiv\bfu_k$, etc.,  passing to the limit $k\to\infty$, and employing \eqref{CoN} along with the trace theorem, we get
$$\ba{cc}\medskip\left.\ba{ll}\medskip
-\Delta \bfv+\nabla p=\lambda\big(\partial_1\bfu-\bfu_0\cdot\nabla\bfu-\bfu\cdot\nabla\bfu_0\big)\\
\Div \bfu=0\ea\right\}\ \ \mbox{in $\Omega$\,,}\\ \medskip
\bfu=\0\ \ \mbox{on $\partial\Omega$\,,}\\
\omega_{\sf n}^2\bfzeta+\varpi\Int{\partial\Omega}{}\mathbb T(\bfv,p)\cdot\bfn=\0\,.
\ea
$$
which proves (b).\hfill$\square$\par\smallskip

Next, we prove

\Bl For any fixed $\lambda>0$ and $\mu\neq0$ the operator
$$
{\mathscr H}_\mu=\mu\,{\sf I}-\lambda\,{\mathscr M}(\lambda)\,
$$
is Fredholm of index 0.
Furthermore, denoting by $\sigma({\mathscr M})$ the spectrum of ${\mathscr M}$, we have that $\sigma({\mathscr M})\cap~ (0,\infty)$ consists at most of a countable number of eigenvalues of finite algebraic multiplicity that can only cluster at 0.
\EL{10}
{\it Proof.}
Let us define $\mathscr T_\mu$ through
 \be{\mathscr H}_\mu=\mu\,{\hat {\mathscr A}}^{-1}\big({\hat {\mathscr A}}+{\frac\lambda\mu}(\hat{\partial_1}+\mathscr P\mathscr C(\lambda)\big):={\hat {\mathscr A}}^{-1}\mathscr T_\mu.\eeq{41}
 Since the operator $\hat{\mathscr A}$ is a homeomorphism by \lemmref{Sto} and ${\mathscr T}_{\mu}$ is Fredholm of index 0 by \lemmref{fredo}, we have
$$
{\rm dim}\,{\sf N}[{\mathscr H}_\mu]={\rm dim}\, {\sf N}[{\mathscr T}_{\mu}]=m<\infty\,.
$$
Moreover, from
$$
\caly={\sf R}({\mathscr T}_\mu)\oplus S_m
$$
with $S_m$  $m$-dimensional subspace, we deduce that for every $\bsfU\in \calz$, we have $\hat{\mathscr A}\bsfU=\bsfU_1+\bsfU_2$, $\bsfU_1\in {\sf R}({\mathscr T}_\mu)$, $\bsfU_2\in S_m$. Therefore, we infer
$$\bsfU={\hat{\mathscr A}}^{-1}\bsfU_1 + {\hat{\mathscr A}}^{-1}\bsfU_2\mbox{, with }\hat{\mathscr A}^{-1}\bsfU_1\in {\sf R}({\mathscr H}_\mu)\mbox{, and }\hat{\mathscr A}^{-1}\bsfU_2\in {\hat{\mathscr A}}^{-1}S_m.$$

It then follows, in particular, that the essential spectrum $\sigma_{\rm ess}({\mathscr M})$ of ${\mathscr M}$ --defined as the set of $\mu$ where ${\mathscr H}_{\mu}$ is not Fredholm-- has empty intersection with $(0,\infty)$.\par
We shall next show that the resolvent set ${\sf P}({\mathscr M})$ of ${\mathscr M}$ has a non-empty intersection with $(0,\infty)$. Since ${\mathscr H}_{\mu}$ is Fredholm of index 0, it is enough to show that, for sufficiently large $\mu>0$, it is ${\sf N}[{\mathscr H}_\mu]=\{0\}$. From \eqref{41}, we see that the latter is equivalent to show that the equation $\mathscr T_\mu(\bsfU)=0$ has only the solution $\bsfU=\0\ni \calx$, for sufficiently large $\mu>0$.
From \eqref{2.2lin}, in turn, this means that the following problem $(\bsfU\equiv(\bfu,\bfchi))$
\be\ba{cc}\medskip\left.\ba{ll}\medskip
-\mu\Delta\bfu+\nabla p=\lambda\,(\partial_1\bfu-\bfu_0\cdot\nabla\bfu-\bfu\cdot\nabla\bfu_0)\\
\Div\bfu=0\ea\right\}\ \ \mbox{in $\Omega$}\,,\\ \medskip
\bfu(x)=\0\,, \ \mbox{ $x\in\partial\Omega$\,,}\\
\omega_{\sf n}^2\bfchi+\varpi\Int{\partial\Omega}{} \mathbb T(\bfu,p)\cdot\bfn=\0\,,
\ea
\eeq{LiPr1}
has only the solution $\bfu=\bfchi=\0$. To show that this is indeed the case, we begin to observe that, since $\bfu_0,\bfu\in L^4(\Omega)$, it follows that  $(\bfu_0\cdot\nabla\bfu+\bfu\cdot\nabla\bfu_0)\in D_0^{-1,2}$. Therefore, from \cite[Theorem VII.7.2]{Gab} we deduce, in particular,
\be
p\in L^2(\Omega)\,.
\eeq{PRS}

Let $\psi_R$ be the cut-off function
introduced in \theoref{noper}. Testing \eqref{LiPr1}$_1$ with $\psi_R\bfu$ and integrating by parts, we get
\be\medskip
\|\psi_R^{\frac12}\nabla{\bfu}\|_2^2=-\half\lambda\,[(\partial_1\psi_R{\bfu},{\bfu})-(|\bfu|^2\bfu_0,\nabla\psi_R)]
-\lambda(\psi_R{\bfu}\cdot\nabla\bfu_0,{\bfu})+({p}\,{\bfu},\nabla\psi_R)\,.
\eeq{FR1}
Thus, passing to the limit $R\to\infty$ in \eqref{FR1}, using \eqref{PRS} and arguing as in the proof of \theoref{noper}, we obtain
$$
\mu\|\nabla\bfu\|_2^2=-\lambda(\bfu_\cdot\nabla\bfu_0,\bfu)\,.
$$
By \theoref{exi}, \eqref{1.9} and H\"older inequality, we deduce
$$
\mu\|\nabla\bfu\|_2^2\le \lambda \|\bfu\|_6^2\|\nabla\bfu_0\|_{\frac32}
\le c\, \lambda \|\nabla\bfu\|_2^2\|\nabla\bfu_0\|_{\frac32}\,.
$$
As a result, if $\mu>c\lambda\|\nabla\bfu_0\|_{\frac32}:=\bar{\mu}$, we conclude $\bfu\equiv\0$, which, by \eqref{LiPr1}$_2$, implies $\bsfU\equiv 0$,
namely, ${\sf P}({\mathscr M})\cap (\bar{\mu},\infty)\neq\emptyset$. Summarizing, we have shown that $\sigma_{\rm ess}({\mathscr M})\cap (0,\infty)=\emptyset$ while ${\sf P}({\mathscr M})\cap (\bar{\mu},\infty)\neq\emptyset$. Therefore, the stated property about eigenvalues is a consequence of classical results in spectral theory \cite[Theorem XVII.2]{GG}.\hfill$\square$\par\smallskip

We now make the following assumption:
\tag{H}
\be
\mbox{the map\ $\lambda\in U(\lambda_s)\mapsto \bfu_0(\lambda)$ is of class $C^2$}\,.
\eeq{H}
This assumption, along with \remref{Ana}, implies that ${\mathscr F}$ is of class $C^2$ in $U\times \mathcal X$. Next,  by \lemmref{10} and \cite[Definition 79.14]{Z1}, for a fixed $\lambda>0$, we  call $\mu\neq 0$  {\it simple eigenvalue} if
\setcounter{equation}{35}
\renewcommand{\theequation}{\arabic{section}.\arabic{subsection}.\arabic{equation}}
\be\ba{ll}\medskip
{\dim}\,{\sf N}[\mu{\sf I}-\lambda\,{\mathscr M}(\lambda)]=1\,;
\\
{\sf N}[\mu{\sf I}-\lambda\,{\mathscr M}(\lambda)]\cap {\sf R}[{\sf I}-\lambda\,{\mathscr M}(\lambda)]=\{0\}\,.
\ea
\eeq{2.2}
As is known, \eqref{2.2}$_2$ can be equivalently reformulated as follows.  Let ${\mathscr M}^*$, be the adjoint of ${\mathscr M}$. Then, from \eqref{2.2}$_1$ and \lemmref{10} we deduce that
$${\rm dim}\,{\sf N}[\mu\,{\sf I}-\lambda\,{\mathscr M}^*(\lambda)]={\rm codim}\,{\sf R}[\mu\,{\sf I}-\lambda\,{\mathscr M}(\lambda)]=1.$$

Indicating by $\bsfW_1\in \calz$ and $\bsfW_1^*\in\calz^{-1}$ two non-zero elements of respectively ${\sf N}[\mu\,{\sf I}-\lambda\,{\mathscr M}(\lambda)]$ and ${\sf N}[\mu\,{\sf I}-\lambda\,{\mathscr M}^*(\lambda)]$, \eqref{2.2}$_2$ is  equivalent (after suitable normalization) to
\be
\langle \bsfW_1^*,\bsfW_1\rangle=1\,,
\eeq{2.3_1111}
where $\langle\cdot,\cdot\rangle$ denotes the duality pairing $\calz\to\calz^{-1}$.\par
The following result holds.
\Bl Suppose  that 1 is a simple eigenvalue of $\lambda_s\,{\mathscr M}(\lambda_s)$ and that \eqref{H} holds. Then, there is  ${U}_0\subseteq U(\lambda_s)$ such that the eigenvalue $\mu=\mu(\lambda)$ of $\lambda\,{\mathscr M}(\lambda)$, $\lambda\in {U}_0$, is simple and of class $C^2$. Moreover,
\be
\mu'(\lambda_s)=-\langle \bsfW_1^*, ({\mathscr M}(\lambda_s)+\lambda_s\,{\mathscr M}^\prime(\lambda_s)\big)(\bsfW_1)\rangle\,,
\eeq{mu1}
where the prime denotes  differentiation with respect to $\lambda$.
\EL{2.1}
{\it Proof.} It can be obtained as a consequence of \eqref{BiP1}--\eqref{BiP2} and \cite[Corollary 79.16]{Z1}.\hfill$\square$\par\smallskip

We are now in a position to prove the main result of this subsection.

\Bt
Suppose that \eqref{H} holds.
If $(\lambda_s,\0)$ is a bifurcation point of \eqref{BiP2}, then the equation
\be
\bsfW-\lambda_s\,{\mathscr M}(\lambda_s)(\bsfW)=\0
\eeq{2.4}
has at least one non-trivial solution $\bsfW_1$. Conversely, assume that $1$ is a simple eigenvalue of $\lambda_s\,{\mathscr M}(\lambda_s)$, namely, \eqref{2.2} holds with $\mu=1$. Then, if $\mu'(\lambda_s)\neq 0$  (transversality condition), there exists exactly one continuous curve of nontrivial solutions to \eqref{BiP1}, $(\bsfU(\lambda),\lambda)\in \calx\times U(\lambda_s)$, with $(\bsfU(\lambda_s),\lambda_s)=(\0,\lambda_s)$.
\ET{2.1.2}
{\it Proof.} Taking into account that \eqref{2.4} is equivalent to ${\rm dim}\,{\sf N}[\mathscr L_s]>0$, the necessary condition follows from \cororef{6.1}; see also \remref{6.}. Moreover, from \eqref{BiP1}--\eqref{BiP2} we have $D_\sub{\tiny \bsfU}{\mathscr F}(\lambda_s,\0)={\sf I}-\lambda_s\mathscr M(\lambda_s)$, which, by \lemmref{10}, is Fredholm of index 0. Therefore, under the assumption $\dim {\sf N}[{\sf I}-\lambda_s\,{\mathscr M}(\lambda_s)]=1$, a classical bifurcation result \cite[Theorem 4.1.12]{Berger} ensures the stated sufficient property provided
$$
D^2_{\lambda\,\mbox{\tiny$\bsfU$}}{\mathscr F}(\lambda_s,\0)\,(\bsfW_1)\not\in {\sf R}[D_\sub{\tiny\bsfU}{\mathscr F}(\lambda_s,\0)]\,,
$$
or, equivalently,
\be
\langle \bsfW_1^*,D^2_{\lambda\,\mbox{\tiny$\bsfU$}}{\mathscr F}(\lambda_s,\0)(\bsfW_1)\rangle\neq 0\,.
\eeq{2.5bis}
By a straightforward computation, from \eqref{BiP1}--\eqref{BiP2} we show that
$$
D^2_{\lambda\,\mbox{\tiny$\bsfU$}}{\mathscr F}(\lambda_s,\0)(\bsfW_1)={\mathscr M}(\lambda_s)(\bsfW_1)+\lambda_s\,{\mathscr M}^\prime(\lambda_s)(\bsfW_1)\,,
$$
so that, if 1 is a simple eigenvalue, by \lemmref{2.1}, condition \eqref{2.5bis} is equivalent to $\mu'(\lambda_s)\neq 0$, which concludes the proof of the theorem.\hfill$\square$\par\smallskip
\Br An equivalent way of expressing \eqref{2.2} is to say that the equation
$$
\bsfW-\lambda_s\mathscr M(\lambda_s)(\bsfW)=\bsfW_1\,,
$$
has no solution. In turn, the latter is equivalent to the condition that the equation
$$
\mathscr L_{\lambda_s}(\bsfW)=\mathscr A(\bsfW_1)
$$
has no solution.
\ER{equiV}

\Br Should the basic flow $\bfu_0(\lambda)$ not depend on $\lambda\in U(\lambda_s)$, then \eqref{mu1} reduces to
$$
\mu'(\lambda_s)=-\langle \bsfW_1^*,{\mathscr M}(\lambda_s)(\bsfW_1)\rangle\,,
$$
which, combined with with  \eqref{2.4} and \eqref{2.3_1111}, furnishes
$$
\mu'(\lambda_s)=-\lambda^{-1}_s\,.
$$
As a result, under the above assumption, the transversality condition is no longer an extra requirement and becomes a consequence of the fact that $\mu=1$ is simple.
\ER{sim}
\Br Also in the light of the previous remarks, we would like to present in more explicit terms the conditions stated in \theoref{2.1.2}, ensuring the occurrence of bifurcation.  Consider the eigenvalue problem
\be\ba{cc}\medskip\left.\ba{ll}\medskip
-\mu(\lambda)\Delta\bfu+\nabla p=\lambda\,(\partial_1\bfu-\bfu_0(\lambda)\cdot\nabla\bfu-\bfu\cdot\nabla\bfu_0(\lambda))\\
\Div\bfu=0\ea\right\}\ \ \mbox{in $\Omega$}\,,\\ \medskip
\bfu(x)=\0\,, \ \mbox{ $x\in\partial\Omega$}\,;\ \
\Lim{|x|\to\infty}\bfu(x)=\0\,,\\
\omega_{\sf n}^2\bfchi+\varpi\Int{\partial\Omega}{} \mathbb T(\bfu,p)\cdot\bfn=\0\,.
\ea
\eeq{LiPr_0}
in the class $(\bfu,\bfchi)\in X^2(\Omega)\times\real^3$. Then,   $(\bfu_0(\lambda_s),\lambda_s)$ is a bifurcation point if the following conditions are met:
\begin{itemize}
\item[(i)] $\mu(\lambda_s)=1$ and the corresponding  eigenspace is one-dimensional, spanned by $(\bfu_1,\bfchi_1)$ ;
\item[(ii)] the problem
$$\ba{cc}\medskip\left.\ba{ll}\medskip
-\Delta\bfu+\nabla p-\lambda_s\,(\partial_1\bfu-\bfu_0(\lambda_s)\cdot\nabla\bfu-\bfu\cdot\nabla\bfu_0(\lambda_s))=\Delta\bfu_1\\
\Div\bfu=0\ea\right\}\ \ \mbox{in $\Omega$}\,,\\ \medskip
\bfu(x)=\0\,, \ \mbox{ $x\in\partial\Omega$}\,;\\
\omega_{\sf n}^2\bfchi+\varpi\Int{\partial\Omega}{} \mathbb T(\bfu,p)\cdot\bfn=\omega_{\sf n}^2\bfchi_1+\varpi\Int{\partial\Omega}{} \mathbb T(\bfu_1,p_1)\cdot\bfn\,,
\ea
$$
has no solution in the class $(\bfu,\bfchi)\in X^2(\Omega)\times\real^3$ ;
\item[(iii)] $\mu'(\lambda)$ satisfies the transversality condition.
\end{itemize}
The last requirement is automatically satisfied  if $\bfu_0(\lambda)\,,\,\lambda\in U(\lambda_s)$, is independent of $\lambda$.
\ER{6F}

\par\medskip\noindent
{\bf Acknowledgments.}
Denis Bonheure is supported by the ARC Advanced 2020-25 ``PDEs in interaction'' and by the Francqui Foundation as Francqui Research Professor 2021-24. Giovanni P. Galdi is  partially supported by the National Science Foundation Grant DMS--2307811. Filippo Gazzola is supported by Dipartimento di Eccellenza 2023-27 of MUR (Italy), PRIN project 2022 ``Partial differential equations and related geometric-functional
inequalities'' (financially supported by the EU, in the framework of the ``Next Generation EU initiative''), and by INdAM.\par\smallskip\noindent
{\bf Declarations.} The authors have no relevant financial or non-financial interests to disclose.
The authors have no competing interests to declare, relevant to the content of this article.
Data sharing not applicable to this article as no datasets were generated or analysed during the current study.

\ed
\begin{thebibliography}{99}

\bibitem{Bab0}Babenko, K.I., On the spectrum of a linearized problem on the flow of a viscous
incompressible fluid around a body (Russian). {\it Dokl. Akad. Nauk SSSR}, {\bf 262}, 64-68 (1982)

\bibitem{BBGGP} Berchio, E., Bonheure, D., Galdi, G.P., Gazzola, F., Perotto, S., Equilibrium configurations of a symmetric body immersed
in a stationary Navier-Stokes flow in a planar channel, preprint (2023)

\bibitem{Berger}Berger, M.S., {\it Nonlinearity and Functional Analysis. Lectures on Nonlinear Problems in
Mathematical Analysis}, Academic Press (1977)

\bibitem{Bev}Blevins, R.D., {\it Flow induced vibrations}, Van Nostrand Reinhold Co., New York (1990)

\bibitem{bocchi} Bocchi, E., Gazzola, F., {\it Asymmetric equilibrium configurations of a body immersed in a 2D laminar flow},
Zeit. Angew. Math. Phys. 74:180 (2023)

\bibitem{BGG}Bonheure, D., Galdi, G.P., Gazzola, F., Equilibrium configuration of a rectangular
obstacle immersed in a channel flow. {\it C. R. Math. Acad. Sci. Paris} {\bf 358}, 887-896 (2020);
updated version in arXiv:2004.10062v2 (2021)

\bibitem{BoGaGaper}Bonheure, D., Galdi, G.P., Gazzola, F., Flow-induced oscillations via Hopf bifurcation in a fluid-solid interaction problem, forthcoming.

\bibitem{BoHiPaSpe} Bonheure, D., Hillairet, M., Patriarca, C., and Sperone, G.,
 Long-time behavior of an anisotropic rigid body interacting with a Poiseuille flow in an unbounded channel, submitted (2023)

\bibitem{Dyr}Dyrbye, C., Hansen, S.O., {\it Wind Loads on Structures}, Wiley, New York (1997)

\bibitem{GaEi}Eiter, T., Galdi, G.P., New results for the Oseen problem with
applications to the Navier-Stokes equations in
exterior domains, RIMS K\^{o}ky\^{u}roku {\bf 2171}, Mathematical Analysis of Viscous Incompressible Fluid (2019) 1-15

\bibitem{FN}Farwig, R., Neustupa, J., Spectral properties in $L^q$ of an Oseen operator modelling
fluid flow past a rotating body, {\it Tohoku Math. J.} {\bf 62}, 287--309 (2010)

\bibitem{Gam} Galdi, G.P., Sulla stabilit\`a incondizionata asintotica in media dei moti stazionari idrodinamici e magnetoidrodinamici in domini limitati o meno, {\it Ricerche Mat.} {\bf 24}, 137-151 (1975)

\bibitem{Gah} Galdi, G.P., On the motion of a rigid body in a viscous liquid: A mathematical analysis with
applications, {\it Handbook of Mathematical Fluid Mechanics}, Elsevier Science, 653-791 (2002)

\bibitem{GaFu} Galdi, G.P., Further properties of steady-state solutions to the
Navier-Stokes problem past a three-dimensional obstacle, J.~Math.~Phys. {\bf 48}, 43 pp. (2007)

\bibitem{Gab}Galdi, G.P., {\it An introduction to the mathematical theory of the Navier-Stokes equations.
Steady-state problems}, Second edition. Springer Monographs in Mathematics, Springer, New York (2011)

\bibitem{GaCe}Galdi, G.P., Steady-state Navier-Stokes problem past a rotating body: geometric-functional properties and related questions., {\it Topics in mathematical fluid mechanics}, 109-197, Springer Lecture Notes in Math., {\bf 2073} (2013)


\bibitem{GaNe} Galdi, G.P., Neustupa, J., Steady-state Navier-Stokes flow around a moving body. {\it Handbook of mathematical analysis in mechanics of viscous fluids}, 341-417, Springer, Cham (2018)

\bibitem{GaSi1}Galdi, G.P., Silvestre, A.L., Strong solutions to the problem of motion of a rigid body in
a Navier-Stokes liquid under the action of prescribed forces and torques, Nonlinear Problems in
Mathematical Physics and Related Topics, I, Int. Math. Ser. (N. Y.), vol. 1, Kluwer/Plenum,
New York, 121-144 (2002)

\bibitem{GaSi2}Galdi, G.P., Silvestre, A.L., On the motion of a rigid body in a Navier-Stokes liquid under the action of a time-periodic force. {\it Indiana Univ. Math. J.} {\bf 58}, 2805-2842 (2009)

\bibitem{GPP}Gazzola, F., Pata, V., Patriarca, C., Attractors for a fluid-structure interaction problem
in a time-dependent phase space, J. Funct. Anal. {\bf 286}, Paper No. 110199, 56 pp. (2024)

\bibitem{GazP}Gazzola, F., Patriarca, C., An explicit threshold for the appearance of lift on the deck of a bridge. {\it J. Math. Fluid Mech.} {\bf 24}, No.1, Paper No. 9, 23 pp. (2022)

\bibitem{GazzSp}Gazzola, F., Sperone, G., Steady Navier-Stokes equations in planar domains with obstacle and explicit bounds for unique solvability, Arch. Ration. Mech. Anal. {\bf 238}, 1283-1347 (2020)

\bibitem{GG}Gohberg, I., Goldberg, S., Kaashoek, M.A., {\it Classes of linear operators: I. Operator Theory},  Advances and Applications, Vol.49, Birkh\"auser Verlag, Basel (1990)

\bibitem{Hey} Heywood, J.G., The Navier-Stokes equations: on the existence, regularity and decay of solutions, {\it Indiana Univ. Math. Journal} \textbf{29}, 639-681 (1980)


\bibitem{PPDL} Paidoussis, M., Price, S., De Langre, E., {\it Fluid-Structure Interactions: Cross-Flow-Induced Instabilities},
Cambridge University Press (2010)

\bibitem{Patri} Patriarca, C., Existence and uniqueness result for a fluid-structure-interaction evolution problem in an unbounded 2D channel, {\it NoDEA Nonlinear Differential Equations Appl.} {\bf 29}, No. 4, Paper No. 39, 38 pp. (2022)


\bibitem{ALS} Silvestre, A.L., On the self-propelled motion of a rigid body in a viscous liquid and on the attainability
of steady symmetric self-propelled motions, {\it J. Math. Fluid Mech.} {\bf 4}, 285-326 (2002)




\bibitem{Will} Williamson, C.H.K., Govardhan, S., Vortex-induced Vibrations, {\it Ann. Rev. Fluid Mech}, {\bf 36}, 413-55 (2004)


\bibitem{Z1}Zeidler, E., {\it Nonlinear Functional Analysis and Applications}, Vol.4, Application to Mathematical Physics, Springer-Verlag, New York (1988)

\end{thebibliography}
